\documentstyle[twoside,11pt]{article}
\title{\huge Spaces of Particles on Manifolds and
Generalized \\ Poincar\'e Dualities
\footnote{This is a revised version of a preprint that
has been circulated in September 96}}

\author{Sadok Kallel\thanks{
The author holds a postdoctoral fellowship with PIms (Pacific
Institute for the mathematical sciences)}\\
{\small\it Pacific Institute for the Mathematical Sciences}\\
{\small\it \& University of British Columbia, Vancouver.}}
\date{}

\begin{document}
\maketitle


\def\la#1{\hbox to #1pc{\leftarrowfill}}
\def\ra#1{\hbox to #1pc{\rightarrowfill}}
\def\fract#1#2{\raise4pt\hbox{$ #1 \atop #2 $}}
\def\decdnar#1{\phantom{\hbox{$\scriptstyle{#1}$}}
\left\downarrow\vbox{\vskip15pt\hbox{$\scriptstyle{#1}$}}\right.}
\def\decupar#1{\phantom{\hbox{$\scriptstyle{#1}$}}
\left\uparrow\vbox{\vskip15pt\hbox{$\scriptstyle{#1}$}}\right.}

\font\sc=cmcsc10 at 11pt

\parskip=1pc
\def\za{\vrule height6pt width4pt depth1pt}

\def\lrar{{\ra 2}}
\def\deg{\hbox{deg}}
\def\loop{\Omega}
\def\cat{\underline{\cal C}}
\def\cont{\underline{\cal S}}
\def\config{C^{\infty}}
\def\party{Par^{\infty}}
\def\sect#1{\hbox{Sec}_{#1}}

\def\sp#1{SP^{#1}}
\def\spy{SP^{\infty}}
\def\tpy{TP^{\infty}}
\def\hol#1{\hbox{Hol}_{#1}}
\def\rat#1{\hbox{Rat}_{#1}}
\def\map#1{\hbox{Map}_{#1}}
\def\div#1{\hbox{Div}^{#1}}

\def\gre{\epsilon}
\def\bbc{{\bf C}}
\def\bbq{{\bf Q}}
\def\bbr{{\bf R}}
\def\bbp{{\bf P}}
\def\bba{{\bf A}}
\def\bbz{{\bf Z}} 
\def\bbf{{\bf F}}


\noindent{\bf\Large\S1 Introduction and Statement of Results}

There are interesting results throughout the literature relating
multi-configuration spaces to mapping spaces (cf. [B], [G], [Gu1-2],
[McD], [S1-2]). In this paper, we use a ``local to global'' {\em
scanning} process based on a construction of Segal to unify
and generalize these results.

First of all by a {\em configuration} on a space $X$ we mean a
collection of unordered points on $X$ (they can be distinct or not).
A multi-configuration will then mean a tuple of configurations with
(possibly) certain relations between them.  Of course, more rigorous
definitions are to follow.

It is known by classical work of G. Segal [S1], that the space of
configurations of {\em distinct} points in Euclidean space is
equivalent in homology to an iterated loop space on a sphere.  Later
work of D. McDuff extended this result to an arbitrary smooth compact
manifold (with boundary) where she showed that the space of
configurations of distinct points there is equivalent in homology to a
space of sections of an appropriate bundle. A bit more later, F. Cohen
and C.F. Bodigheimer proved a similar result for spaces of
configurations of distinct points {\em with labels} (see [B]).

Both Segal and McDuff extended their ideas to spaces made out of pairs
of configurations. While Segal worked with {\em divisor spaces} made
out of pairs of configurations having no points in common on a
punctured Riemann surface [S2], McDuff dealt with what she coined the
{\em space of positive and negative particles} on a general smooth
manifold. Both were able to identify these spaces with some function
spaces.

This paper extends and generalizes the work of Segal and McDuff in a
great many directions. It also sets a context in which these types of
results can be viewed and interpreted by relating them to more
classical aspects of algebraic topology, as well as to some recent
problems stemming from Gauge theory and dealing with the topology of
holomorphic mapping spaces.

A starting point for us has been to address the following question:
which (multi-) configuration spaces can be used to model mapping
spaces (and vice-versa). Such considerations have led us to introduce
a new and general class of multi-configurations, the {\em particle
spaces} and these include symmetric and truncated products, divisor
spaces, configurations of distinct points, spaces of positive and
negative particles, and most other known examples in the literature.

Our basic definition is: {\sl A particle space is a
multi-configuration space with a partial monoid structure.}  If one
defines the support of a multi-configuration to mean the locus of the
points making up the multi-configuration, then the partial monoid
structure will be ``concatenation'' defined on multiconfigurations
having disjoint support.

Given a manifold $M$, the most basic example of a particle space on
$M$ is the infinite symmetric product $\spy (M)=\coprod\sp{n}(M)$ (and
this is a monoid).  Another standard particle space is the
(traditional) configuration space $C^{\infty}(M)\subset\spy (M)$
consisting of unordered disjoint points of $M$ (and inheriting a
partial monoid structure). We agree on the following notation: an
element $\zeta$ in $\sp{n}(M)$ can be written both as the formal sum
$\sum n_ix_i, x_i\in M$, $n_i\in {\bf N}$ and $\sum n_i=n$, or as an
unordered tuple $\langle x_1,\ldots, x_n\rangle$.

Generally we define $\party (M)$ to be any quotient of any subset of a
finite product $\prod\spy (M)$ satisfying the partial monoid structure
referred to above. A point of $\party (M)$ will then be a
multiconfiguration on $M$ with certain constraints. The particle
spaces are defined on any manifold $M$ (which we do assume in this
paper to be smooth) and hence we can talk of a {\sl particle functor}
$\party$.  Here are some examples of those functors and spaces we
study in this paper:

\noindent $\bullet$
Symmetric product spaces with ``bounded multiplicity''. Given $M$ as
above and an integer $d\geq 1$, we define
$$\spy_d(M) =\{\sum n_ix_i\in\spy (M)~|~n_i\leq d\}$$
Of course $\spy_1 = \config$ is the configuration space of distinct
points. \hfill\break
\noindent $\bullet$
$\party (M) = \{(\zeta_1,\ldots,\zeta_k)\in\spy
(M)^k~|~ \zeta_i\cap\zeta_j=\emptyset, i\neq j\}.$ 
A related space
will be the set of $k$-tuples of configurations $j$ of which are distinct,
$j\leq k$. \hfill\break
\noindent$\bullet$ $\party (M) = 
\prod^n\spy (M)/\Delta (\spy (M))$
where $\Delta$ is the submonoid generated by diagonal elements.
\hfill\break
\noindent$\bullet$ Truncated symmetric products and these refer to
$\tpy_p(M)=\spy (M)/x^p\sim *$ (here we're thinking of $\spy (M)$ as a
topological monoid with $*\in M$ the identity element).\hfill\break 
\noindent$\bullet$  Spaces of positive and negative particles of
McDuff and these refer to
$\party (M) = C^{+\over }(M)=C(M)\times C(M)/\sim$ where $\sim$
is the identification
$$
(\zeta_1,\eta_1)\sim_R(\zeta_2,\eta_2)\Leftrightarrow \zeta_1-\eta_1
=\zeta_2-\eta_2.
$$ 
\noindent$\bullet$ The divisor spaces of Segal studied in connection
with the space of holomorphic maps of Riemann surfaces into projective
spaces (see [K1] and [C$^2$M$^2$]). They are defined as
$$\div{n}(M)=\{(\zeta_1,\ldots,
\zeta_n)\in\spy (M)^n~ |~\zeta_1\cap\zeta_2\cap\cdots
\cap\zeta_n=\emptyset\}.
$$ 

The key property of the particle spaces is that when you look closely
at a multiconfiguration of $\party (M)$ in the neighborhood $D$ of a
point (that is when you {\sl scan} the manifold), what you see
is a multiconfiguration living in $\party (D)$. This restriction
property turns out to be a direct consequence of the partial monoid
structure put on $\party (M)$.

As is standard, one can define relative particle spaces whereby the
functor $\party$ can be applied to a pair of spaces. If $N\subset M$,
then $\party (M,N)$ consists (roughly) of all those multiconfigurations
in $\party (M-N)$ which get identified 
as they approach $N$. It is not hard to see that the scanning
property mentioned in the previous paragraph establishes (at least
for parallelizable manifolds $M$) the existence of a map
$$S: \party (M)\lrar \map{}(M, \party(S^n,*))$$
where $*\in S^n$ can be chosen to be the north pole.

For a given space $M$, $\party (M)$ is a disconnected partial monoid
(with components not very comparable.) It turns out that by ``group
completing'' with respect to this partial monoid structure, one
obtains a space $Par (M)$ which is better behaved (and all of whose
components are homeomorphic). The functor $Par$ (which we construct in
\S7) is the last ingredient we need and we are now in a position to
state the main result of this paper.

\noindent{\bf Main Theorem} 1.1:~{\sl Let $M$ be an $n$ dimensional,
smooth, compact (possibly with boundary) and connected manifold. 
Then there is a fibre bundle
$$\party (S^n,*)\lrar E_{\party}\lrar M\leqno{1.2}$$
with a (zero) section. Choose
$N$ to be a closed ANR in $M$ and assume that either $N\neq\emptyset$
or $\partial M\neq\emptyset$. 
Then there is a homology equivalence (induced by scanning)
$$
S_*: H_*\left(Par (M-N); \bbz \right) \fract{\cong}{\ra 2}
H_*\left(\sect{}(M,N\cup\partial M, \party (S^n,*));\bbz\right)
$$
where $\sect{} (M,A,\party (S^n,*))$ is the space of sections
of $1.2$ trivial over $A$.}

The above theorem has several variants described throughout this
paper. An immediate question one asks is when can the homology
equivalence of theorem 1.1 be upgraded to a homotopy equivalence. 
We resolve this as follows.

\noindent{\bf Theorem} 1.3:~{\sl Let $N,M$ be as in 1.1 and suppose
$\pi_1(Par(\bbr^n))$ is abelian, then scanning is a homotopy equivalence
$$Par (M-N) \fract{\simeq}{\ra 2} \sect{}(M,N\cup\partial M, \party
(S^n,*)).$$}

\noindent{\bf Corollaries and Examples}: \hfill\break 
$\bullet$ When $M$ is parallelizable, the bundle of
configurations $1.2$ trivializes and sections turn into
maps into the fiber. One therefore has the equivalence
$$ H_*\left(Par (M-N); \bbz \right) \fract{\cong}{\ra 2}
H_*\left(\map{}(M, N,\party (S^n,*));\bbz\right)$$
where $\map{}(M,N,\party (S^n,*))$ is the space of (based) maps
sending $N$ to the canonical basepoint in $\party (S^n,*)$.
When $N=*$, we write $\map{}^*(M,\party (S^n,*))$ for the corresponding
mapping space.
\hfill\break
$\bullet$~(Segal [S1])~Let $M_g$ be a genus $g$ Riemann
surface. Then $\div{2}(M_g-*)\simeq\map{c}^*(M,\bbp\vee\bbp )$ where
$\bbp=K(\bbz,2)$ is the infinite complex projective space and where
$\map{c}^*$ is any component of the subspace of based maps (see 11.4).
\hfill\break
$\bullet$~(McDuff [McD1])~
$C^{+\over}(\bbr^n)\simeq\Omega^n\left(S^n\times
S^n/\Delta\right)$ where $\Delta$ is the diagonal copy of $S^n$ in
$S^n\times S^n$.\hfill\break
$\bullet$ Let $C$ be the configuration functor associated to $\config$.
Then
$$H_*(C(\bbr^n);\bbz )\fract{\cong}{\ra 2} H_*(\Omega^n(S^n);\bbz )
~~~~~\hbox{[S2]}\leqno{1.4}.$$
$\bullet$ We can generalize 1.4 as follows. Let
$C^{(k)}(\bbr^n)\subset\prod^kC (\bbr^n)$
consist of the subspace of {\em pairwise disjoint} configurations. Then
$$H_*(C^{(k)}(\bbr^n);\bbz )\fract{\cong}{\ra 2}
H_*(\Omega^n (\underbrace{S^n\vee\cdots\vee S^n}_k);\bbz).$$

Other interesting examples we discuss are the symmetric products with
bounded multiplicity which we introduced earlier and denoted by
$\spy_d(M), d\geq 1$. For $M$ either open or with boundary,
the ``partial'' completion $SP_d (M)$ has the following very simple
description. Choose an end (or a tubular neighborhood of the boundary)
and construct a nested sequence $\{U_i\}$ of neighborhoods of it. By
choosing a sequence of disjoint points $z_i\in U_i-U_{i+1}$ we obtain
maps
$$\spy_d(M-U_i)\fract{+z_i}{\ra 2}\spy_d(M-U_{i+1})\fract{+z_{i+1}}{\ra 2}
\cdots$$
and the direct limit we denote by $SP_d(M)$. The following
proposition, which we state in the special case $M=\bbr^n$ (see
theorem 11.7), is a direct corollary of 1.1 and 1.3 once we observe
that $\spy_d(S^n,*)\simeq\sp{d}(S^n)$ (cf. \S11).

\noindent{\bf Proposition} 1.5:~{\sl Scanning $S$ is a homotopy equivalence
$$SP_d (\bbr^n)\fract{\simeq}{\ra 3} \Omega^n_0\sp{d}(S^n)$$
whenever $d>1$, and a homology equivalence when $d=1$.}

\noindent{\bf Note}: The proposition above has also been obtained by
M. Guest, A. Kozlowski and K. Yamaguchi [GKY] (who state it for the
case $n=2$; cf.  \S11.2).  A labelled analog of it is given in [K2]
and yields a direct generalization of the May-Milgram model for
iterated loop spaces.  One might note that 1.5 provides yet another
extension of Segal's result (1.4).

One main interest in theorem 1.1 is the way it relates to and
generalizes many of the classical dualities on manifolds. The
following theorem (obtained earlier by Pawel Gajer [G] using different
techniques) is obtained after a close analysis of the bundle 1.2 for
the case $\party =\spy$.

\noindent{\bf Theorem} 1.6:~{\sl Let $M$ be $n$ dimensional, smooth
and compact, and let $N$ be an ANR in $M$. Suppose $M$ is orientable.
Then scanning induces a homotopy equivalence
$$
S:\spy (M-N,*)\fract{\simeq}{\ra 2}\map{c}(M, N\cup\partial M,\spy (S^n,*))
$$
where $\map{c}$ is any component of the space of (based) maps}.

\noindent{\bf Corollary} 1.7 (Alexander-Lefshetz-Poincar\'e):~{\sl 
Let $M$ and $N$ be as above, then
$${\tilde H}_*(M-N;\bbz )\cong H^{n-*}(M, N\cup\partial M,\bbz ).$$}

This work finds its origins in an attempt to construct {\em
configuration space models} for spaces of holomorphic maps on Riemann
surfaces $M_g$. In the past decade and as a result of the increasing
``rapprochement'' between mathematics and physics, there has been a
flurry of activity towards understanding the topology of spaces
$\hol{}^*(M_g,X)$ of (based) holomorphic functions into various
algebraic varieties. The general picture that emerges there is that
for many special rational $X$'s one has the following relationships
(eg. [Gu1-2], [BHMM], [KM])

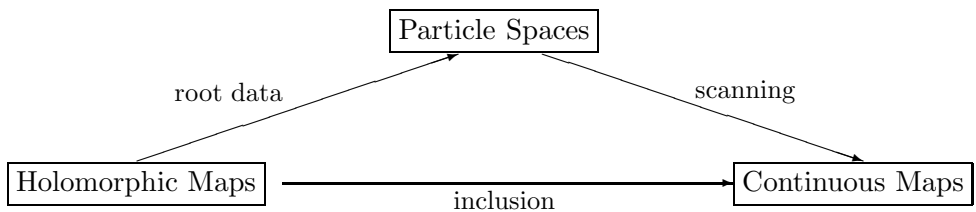
\setlength{\unitlength}{.4 truein}
\begin{picture}(6,3.6)(0,0.3)
\put(6,3){\fbox{Particle Spaces}}
\put(1,1){\fbox{Holomorphic Maps}}
\put(10.5,1){\fbox{Continuous Maps}}
\put(4.6,1.1){\vector(1,0){5.9}}
\put(2.7,1.4){\vector(3,1){4.2}}
\put(8,2.8){\vector(3,-1){4.2}}
\put(10,2.3){\makebox(0,0)[l]{\small scanning}}
\put(7.5,1){\makebox(0,0)[t]{\small inclusion}}
\put(4.6,2.3){\makebox(0,0)[r]{\small root data}}
\end{picture}

In this framework, one uses the particle spaces to provide models for
spaces of holomorphic maps on a Riemann surface, which themselves are
{\em suitable approximations} to spaces of all maps. The
following (perhaps unsuspected) corollary is given in \S14 and it
recuperates a known theorem of M. Guest [Gu2].

\noindent{\bf Corollary} 1.8:~(Guest)~{\sl Let $X$ be a projective
toric variety (non-singular). The natural inclusions
$i_D:\hol{D}^*(S^2,V)\lrar\loop^2_DV$ (where $D$ are multidegrees
depending on $V$) induce a homotopy equivalence when
$D$ goes to $\infty$; i.e.  
$$\displaystyle \lim_{D\rightarrow\infty}\hol{D}^*(S^2,V)
\fract{\simeq}{\lrar}\loop^2_0V$$
where $\loop^2_0V$ is any component of $\loop^2V$.}

\noindent{\sc Remark} 1.9: The equivalence above between spaces of
rational maps and loop spaces of certain projective varieties has been
observed initially by Segal for the case of $V=\bbp^n$ and later
extended to more general flag manifolds by several authors (see
[C$^2$M$^2$], [BHMM] and references therein). In light of the methods
used in this paper, it turns out that it is precisely the partial
monoid structure that is exhibited by the root data of rational maps
on toric varieties that induces the equivalence with the second loop
space of $V$.  This shouldn't be surprising in light of earlier work
of Segal ([S4]) and should provide an interesting insight into why
equivalences of the sort should hold.

Finally, it is not hard to see that the ideas presented above apply
equally well (but in a different context) to obtaining space level
descriptions of Spanier-Whitehead duality for any generalized homology
theory (cf. \S15).

\noindent{\bf Theorem} 1.10:~{\sl Let $E$ be a connected $\Omega$
spectrum and define the functor $F_{\bf E}(-)=\Omega^{\infty}(\bf
E\wedge -)$ on the category of CW complexes.  Then for all $X\in
CW$, there is a homotopy equivalence
$$
S:F(X)\fract{\simeq}{\ra 2}\map{*}(D(X,k), F(S^k))
$$
where $D(X,k)=S^k-X$ is the Spanier-Whitehead dual of $X\hookrightarrow
S^k$.}

\noindent{\bf Corollary} 1.11:~(Spanier-Whitehead duality)~ {\sl Let
$h$ be any homology theory and suppose $A,B\in S^k$, $A$ and $B$ are
$n$ dual. Then there is an isomorphism 
$$h_i(B)\cong h^{n-1-i}(A).$$}

\noindent{\sc Acknowledgements:} The author is grateful to Pawel
Gajer for commenting on an early version of this paper.  He is also
grateful to M. Guest, A. Kozlowski and K. Yamaguchi for sending him
their preprint [GKY]. We thank the Fields institute for its
hospitality at the time this work was being conducted.  The author
would like to express his gratitude to the organizers of the special
homotopy theory program for the 95/96 theme year at the Fields
institute; in particular to R. Kane and P. Selick.  Finally, he wishes
to acknowledge the generous support of both the CRM in Montr\'eal and
of Professor J. Hurtubise (McGill).


\vskip 10pt
\noindent{\bf\Large\S2 Quasifibrations and Homology Fibrations}

\noindent{\sc Definition}: 
Recall a map $f: Y\rightarrow X$ is a quasifibration if $\forall x\in
X$ the inclusion of $f^{-1}(x)$ into the homotopy fiber over $x$ is a
weak homotopy equivalence.  Roughly speaking, while a fibration enjoys
the property of homotopy lifting ``on the nose'', one may need to
``deform'' homotopies before being able to lift them for the case of a
quasifibration.  A standard example is given by the projection $\pi$
depicted below

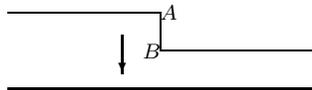
\setlength{\unitlength}{.4 truein}
\begin{picture}(6,2.5)(-4,0)
\put(0,2){\line(1,0){2}}
\put(2,2){\line(0,-1){0.5}}
\put(2,1.5){\line(1,0){2}}
\put(0,1){\line(1,0){4}}
\put(2,2){\makebox(0,0)[l]{$\scriptstyle A$}}
\put(2,1.5){\makebox(0,0)[r]{$\scriptstyle B$}}
\put(1.5,1.7){\vector(0,-1){.5}}
\end{picture}

Clearly, $\pi$ is not a fibration for one cannot lift homotopies that
don't ``spend much time" at the point $\{{1\over 2}\}=\pi ([AB])$.
This projection is however a quasifibration and by allowing homotopies
to ``live a while" over certain closed sets ($\{{1\over 2}\}$ in this
case) one should be able to lift them. This is the essense of 2.2
below.

\noindent{\sc Definition} 2.1:~{\sl
A map $\pi:E\rightarrow B$ is a homology fibration
if for each $b\in B$, the natural map $\pi^{-1}(b)\rightarrow F$ into
the homotopy fiber is a homology equivalence.}

The general criterion developed by Dold and Thom in [DT]
to show that a map is a quasifibration can be extended to 
include the case of homology fibrations as well. This gives

\noindent{\bf Criterion} 2.2: [DT]~{\sl Let $X_1\subset\cdots\subset
X_k\subset\cdots$ be a (finite) filtration of $X$ by closed subspaces and
let $f:Y\lrar X$ be a map satisfying\hfil\break
\hbox to 0.3in{\hfill} (i) $f$ is a fibre bundle over $X_{k+1}-X_k$
with fiber $F$\hfil\break
\hbox to 0.3in{\hfill} (ii) There is an open set $X_{k-1}\subset
U_k\subset X_k$ and a deformation retraction of $r_t$ of
$U_k$ to $X_{k-1}$ which can be lifted to a deformation retraction 
${\tilde r}_t$ (upstairs) of $f^{-1}(U_k)$ to $f^{-1}(X_{k-1})$.\hfil\break
\hbox to 0.2in{\hfill}
The map $f$ is a quasifibration (resp. homology fibration)
if ${\tilde r}_1:f^{-1}(x)\rightarrow f^{-1}(r_1(x))$ is a weak
homotopy (resp. homology)
equivalence for all $x\in U_k$.}

\noindent{\sc Terminology}: The maps ${\tilde r}_1$ are referred to by
McDuff as {\em attaching maps}. We will adopt the same terminology.

\noindent{\bf Remark} 2.3: A slightly more general version of 2.2
holds: Suppose $X_{\vec k} = X_{k_1,\ldots, k_n}$ is a cover of $X$ by
closed subsets such that $X_{k_1,\ldots, k_n}\subset X_{k_1,\ldots,
k_{i+1},\ldots, k_n}$ for all $1\leq i\leq n$ and suppose that $\pi$
is a trivial fibration over $X_{k_1,\ldots, k_n}-\bigcup_i
X_{k_1,\ldots, k_{i-1},\ldots , k_n}$.  Let $U_{\vec k}$ and the
attaching maps ${\tilde r}_1$ be defined as in 2.2 where $U_{\vec k}$
retracts down to $X_{\vec k}$ via a retraction $r$.  
The very same criterion as in 2.2 states that
if the attaching maps over the $X_{\vec k}$ are homotopy
(resp. homology) equivalences, then $\pi$ is a quasifibration (resp. a
homology fibration).

\noindent{\sc Example} 2.4: We illustrate 2.2 for the projection
$\pi:Y\rightarrow [0,1]$ depicted earlier.  Set $X_1={1\over 2}$ and
$X_2=[0,1]$. One can let $U_1$ be $(0,1)$. The retraction $r$ over the
time interval $[0,2]$ is chosen to shrink linearly $(0,1)$ to
$\{{1\over 2}\}$ over the time interval $[0,1]$ and to be stationary
at ${1\over 2}$ for $t\in [1,2]$.  Now ${\tilde r}_1$ corresponds to
the following: it shrinks $Y$ to the vertical line segment $AB$
(linearly over the time interval $[0,1]$) leaving $AB$ fixed. Then for
$t\in [1,2]$, it slides the point $A$ to the end point $B$ along
$AB$. By construction, this gives a lift $\tilde r$ of $r$ and we have
that ${\tilde r}_1(r^{-1}(U_1))=\{B\}$.

\noindent{\sc Example} 2.5: We defined earlier $\spy (M) =
\coprod\sp{n}(M)$.  Choose a basepoint $*\in M$ and construct
inclusions $\sp{n}(M)\hookrightarrow \sp{n+1}(M)$ given by adjoining
basepoint $\sum n_ix_i\mapsto \sim n_ix_i+*$.  The direct limit is
denoted by $\spy (M,*)$. The following classical theorem of Dold-Thom
will serve as a prototype for later proofs.

\noindent{\sc Definition} 2.6: For $*\in N\subset M$, we define
$\spy (M,N)$ to be $\spy (M/N,*)$. An equivalent description of this
space is given in \S6.

\noindent{\bf Proposition} 2.7:~{\sl Let $N\hookrightarrow M\rightarrow M/N$
be a cofibration and choose a basepoint $*\in N\subset M$. Then
$$\spy (N,*)\lrar\spy (M,*)\lrar\spy (M,N)$$
is a quasifibration.}

\noindent{\sc Proof:} Let $X^k = \sp{k}(M,N)$ be the image of $\sp{k}(M)$
under the quotient map $\sp{k}(M)\rightarrow\sp{k}(M/N)\hookrightarrow
\spy (M,N)$. It should be clear that
$$X_k=\{D\in\spy (M,N)~|~card(D\cap (M-N))\leq k\}$$
and that the $X_k$ provide an increasing filtration of the $\spy (M,N)$.
Since $N\hookrightarrow M$ is a cofibration, there is a neighborhood
retract $U$ of $N$ in $M$; that is there is an open $N\subset U\subset M$
and a continuous $r:M\lrar M$ such that $r$ leaves $M-U$ and $N$ invariant and
maps $U$ to $N$. The map $r$ lifts (additively) to $\spy(M)$ and we write 
this map as $\tilde r$. Let now
$$U_k=\{D\in\spy (M,N)~|~card (D\cap (M-N))\leq k~\hbox{and at least one
element of $D$ is in $U$}\}.$$
Clearly $X_{k-1}\subset U_k\subset X_k$ and $r_{|U_k}$ retracts $U_k$
onto $X_{k-1}$. It is also clear that over 
$$X_k-X_{k-1}=\{D\in\spy (M,N)~|~card (D\cap (M-N))=k\}$$
the projection $\pi:\spy (M)\rightarrow\spy (M,N)$ is a product
$\pi^{-1}(X_k-X_{k-1})=\spy (N)\times (X_k-X_{k-1})$ and hence is 
trivial there. It remains then to check condition (ii) of 2.2.

Let $b=\langle z_1,\ldots, z_k\rangle\in X_k$ and $z_i\in M-N$. 
Then $b\in U_k$ if one of the
$z_i$ is in $U$. Write $b=\langle z_1,\ldots, z_1',\ldots\rangle$ where
$z_i\in M-U$ and $z_i'\in U-N$. Let $F$ denote $\spy (N)$. One uses
the trivialization of $\pi$ over $X_k-X_{k-1}$ to write
$\pi^{-1}(b)=b+F$ and so the lifted retraction ${\tilde r}:\pi^{-1}(b)
\lrar\pi^{-1}(r(b))$ takes the form
\begin{eqnarray*}
   {\tilde r}_1: \langle z_1,\ldots, z_1',\ldots \rangle + F&\lrar&
   \langle r(z_1),\ldots, r(z_1'),\ldots \rangle + F\\
   &\lrar&\langle r(z_1),\ldots\rangle + (\langle r(z_1'),\ldots\rangle + F)
\end{eqnarray*}
But since $r_1(z_i')\in N,\forall i$, they can be connected to basepoint by paths
and this defines a homotopy of $(\langle r(z_1'),\ldots\rangle + F)\simeq
F$ and hence ${\tilde r}_1(\pi^{-1}(b)\simeq \langle r_1(z_1),\ldots\rangle + F
\simeq\pi^{-1}(r_1(b))$.
\hfill\za


\vskip 10pt
\noindent{\bf\Large\S3 Particle Functors and Particle Spaces}

In this section, we define the $\party$ spaces associated to path
connected spaces $M$.  The starting point here is the
multiconfiguration space $\prod^k\spy (M)$.  We are interested in
subsets and quotients of this space satisfying an ``adjunction"
condition. 

\noindent{\sc Terminology:} 
Consider the $n$-tuple of configurations ${\vec\zeta} =
(\zeta_1,\ldots,\zeta_k)\in\spy (M)^{\times n}$. \hfil\break
\noindent $\bullet$ 
The support of $\zeta_i$ is the set of points making up $\zeta_i$
and the support of $\zeta$ is the union of the supports of the $\zeta_i$.
\hfill\break
\noindent $\bullet$
A {\em subtuple} $\vec\zeta '$ of $\vec\zeta$ consists of a $k$-tuple
$(\zeta'_1,\ldots ,\zeta'_k)$ of subconfigurations
$\zeta'_i\subset\zeta_i$.\hfil\break
\noindent $\bullet$ 
${\vec\zeta}$ is said to lie 
in $A\subset M$ if the support of $\zeta$ is in $A$.
Equivalently ${\vec\zeta}\in\spy (A)^k\subset\spy (M)^k$.  \hfil\break
\noindent $\bullet$ 
Two configurations $\vec\zeta$ and $\vec\zeta'$ are {\em distinct} 
if their supports are 
distincts.  They are {\em disjoint} if they lie in disjoint subsets of $M$;
i.e. if $\vec\zeta\cap \vec\zeta'=\emptyset$. Of course disjoint implies
distinct.\hfil\break
\noindent $\bullet$ 
Given $A\subset\spy (M)^k$ and ${\vec\zeta}\in\spy (M)^k$ then
${\vec\zeta}\cap A$ is the subtuple of $\vec\zeta$ made out of
the points of $\vec\zeta$ that are in $A$.

\noindent{\sc Definition and Notation}: $\prod\spy (M)=\spy (M)^k$
is a topological monoid and we write its pairing as $+$.  We denote by
$\cat$ be the category of spaces with {\em injections} as morphisms.

\noindent{\sc Definition}~3.1 ({\bf Particle Spaces of the first
kind}):~ A (sub) {\em particle} functor $\party$, or a particle space
of the {\sl first kind}, is a {\em covariant} functor $\cat\rightarrow
\cat$ satisfying the following two properties:
$$
\party: M\mapsto \party (M)\subset\spy (M)^k,~
\hbox{for some}~k>0, \forall M\in\cat\leqno{{\bf P}1}
$$
and $\forall A,B\subset M\in\cat, A\cap B=\emptyset$, the symmetric product
pairing $+$ yields an identification
$$\party (A\sqcup B) = \party (A)+\party (B) \subset \party (M).
\leqno{{\bf P}2}$$

\noindent{\bf Remark} 3.2: The symmetric product pairing $+$ restricts
to a partial pairing on $\party$ whereby element can be added only if
they have disjoint support.  This endows $\party (M)$ with a {\em
partial monoid} structure.  Note that the functoriality of $\party$
implies the following naturality for all $N\subset M$
$$\matrix{
   \party (N)&\hookrightarrow&\party (M)\cr
   \decdnar{\subset}&&\decdnar{\subset}\cr
   \spy (N)^k&\hookrightarrow&\spy (M)^k.
\cr}$$

\noindent{\bf Example} 3.3: Let $A\subset M\in\cat$, $A$ open, and let
$F(M)$ be the space
$$
F(M)=\left\{(\zeta_1,\zeta_2,\zeta_3)\in\spy (M)^3~|~\zeta_1\cap
A=\emptyset.  \right\}
$$
$F(M)$ is not a particle space since it is not induced from a 
functor. 

\noindent{\bf Example} 3.4: Consider the space $F(M)\subset \spy
(M)^3$ consisting of triples $(\zeta_1,\zeta_2,\zeta_3)$ such that
$$\deg (\zeta_1) = \deg (\zeta_2) +\deg (\zeta_3).$$ Then $F$ defines
a functor $\cat\rightarrow\cat$. It however doesn't satisfy {\bf P}2
for it is easy to see that the inclusion $F(A)+F(B)\subset F(A\sqcup
B)$ is proper.

The following gives a description of particle functors of the first kind.

\noindent{\bf Lemma} 3.5:~{\sl Let $F$ be a functor
$\cat\rightarrow\cat$ satisfying {\bf P}1; that is $\exists k>0$ such
that $F(M)\subset\spy (M)^k$ for all $M\in\cat$.  Then $F$ is a
$\party$ functor if and only if for all $N\subset M$
$$F (N) = F (M)\cap\spy (N)^k.$$}

\noindent{\sc Proof:} Let $A$ and $B$ be disjoint in $M$ and suppose
$F(A)=F(M)\cap\spy (A)^k$ (same for $B$). Then
\begin{eqnarray*}
    F(A\sqcup B)&=&F(A\sqcup B)\cap\spy (A\sqcup B)^k\\
    &=&F(A\sqcup B)\cap (\spy (A)^k\times\spy (B)^k)=
    F(A\sqcup B)\cap (\spy (A)^k+\spy (B)^k)\\
    &=&F(A\sqcup B)\cap\spy (A)^k + F(A\sqcup B)\cap\spy (B)^k\\
    &=&F(A)+F(B)
\end{eqnarray*}
and $F$ is indeed a $\party$ functor. Suppose now that $F=\party$ for
some particle functor and let $N\subset M\in\cat$. Then
\begin{eqnarray*}
   \party(M)\cap\spy (N)^k& =& (\party(N)+\party(M-N))\cap\spy (N)^k\\
   &=&\party(N)\cap\spy (N)^k + \party(M-N)\cap\spy (N)^k\\
   &=&\party(N)\cap\spy (N)^k=\party(N)
\end{eqnarray*}
as desired and this proves the lemma.
\hfill\za

\noindent{\sc Definition} 3.6 ({\bf Particle functors}):~ A functor
$\party:\cat\rightarrow\cat$ is a particle functor if there are
quotient maps
$$q_M:\party_1(M)\lrar\party (M),~\forall M\in\cat\leqno{{\bf Q}1}$$
of some particle space of the first kind such that the partial pairing
$+$ on $\party_1$ descends to a partial pairing of quotient spaces; i.e.
$$\matrix{
   \party_1 (A)\times\party_1 (B)&\fract{+}{\ra 2}&\party_1 (M)\cr
   \decdnar{q_A\times q_B}&&\decdnar{q_M}\cr
   \party (A)\times\party (B)&\fract{+}{\ra 2}&\party (M)
\cr}$$
whenever $A,B\subset M$, $A\cap B=\emptyset$.

\noindent{\bf Remark} 3.7: When $\party$ is a quotient of $\prod\spy$,
then $A$ and $B$ don't need to be distinct and we can demand that $+$
commutes with $q$ i.e
$$\matrix{
   \spy (M)^k\times\spy (M)^k&\fract{+}{\ra 2}&\spy (M)^k\cr
   \decdnar{q\times q}&&\decdnar{q}\cr
   \party (M)\times\party (M)&\fract{+}{\ra 2}&\party (M)
\cr}$$
and so $\party (M)$ in this case has automatically a monoidal
structure given by $+$ above.

\noindent{\sc Notation}: We write an element $\zeta\in\party (M)$ as
a tuple $(\zeta_1,\ldots, \zeta_n)$ which could either be in
$\spy (M)^k$ or could represent $q^{-1}(\zeta)$ in $\spy (M)^k$.

\noindent{\bf Remark} 3.8: It has been pointed out to the author that
definition 3.6 is closely related to a similar definition of M. Weiss
dealing with spaces of immersions and embeddings of manifolds
(eg. preprint ``Embeddings from the point of view of Immersion
theory''). Unfortunately we are not very knowledgeable of the work of
Weiss at this point to make the analogy precise.


\vskip 10pt
\noindent{\bf\Large\S4 Construction of Particle Spaces}

\noindent{\bf Definition} 4.1:~ Let $U$ be a topological (partial)
monoid and $A\subset U$ any subspace. Then by $U//A$ we mean the
identification space
$$U//A= U/a+x\sim x,~a\in A~\hbox{and whenever}~a+x~\hbox{is defined}.$$
If $U$ is a monoid and $A$ a submonoid, then $U//A$ is simply the quotient
monoid.

\noindent{\bf Lemma} 4.2:~{\sl Let $\party_1$ and $\party_2$ be two
particle functors, $M\in\cat$. \hfill\break
(a) Then $\party_1(M)\times\party_2(M)$ is a particle space and
$\party_1\times\party_2$ a particle functor. \hfill\break
(b) If $\party_1(M)\subset\party_2(M)$, then
$\party_1 (M)//\party_2 (M)$ is a particle space.}

\noindent{\sc Proof:} We observe that
\begin{eqnarray*}
   \party_1 (A\sqcup B)^k//\party_2 (A\sqcup B)
   &=&[\party_1 (A)^k+\party_1 (B)^k]//[\party_2 (A)+\party_2 (B)]\\
   &=&\party_1 (A)^k//\party_2 (A)+\party_1 (B)//\party_2 (B)
\end{eqnarray*}
and (b) follows.
\hfill\za

\noindent{\sc Example} 4.3:~A space we looked at before is
$C^{+\over}$ which is given as the quotient of
$\config(M)^2\subset\spy (M)^2$ by $\party_2(M)=\Delta\config (M)$
where $\Delta$ is the diagonal
$$\Delta:\config (M)\lrar\config (M)\times\config (M),~\zeta
\mapsto (\zeta,\zeta).$$

Consider now for each $M\in\cat$ a map of monoids $f_M:\spy
(M)^m\lrar\spy (M)^n$, $m, n$ positive integers. We assume that
the maps $f_M, M\in\cat$ are compatible with inclusions
$N\subset M$; that is there are commutative diagrams
$$
\begin{array}{ccc}
   \spy (N)^{m}&\fract{f_N}{\ra 2}&\spy (N)^{n}\\
   \decdnar{\subset}&&\decdnar{\subset}\\
   \spy (M)^{m}&\fract{f_M}{\ra 2}&\spy (M)^{n}.
   \end{array}\leqno{4.4}
$$
\noindent{\bf Definition} 4.5: Given a subset $\emptyset\neq A\subset
\prod^n\spy (M)=\spy (M)^{\times n}$ we denote by $(A)\in\spy (M)$ the
submonoid
$$(A) = \{ a+x, a\in A, x\in\spy (M)^{\times n}\}.$$

\noindent{\bf Proposition} 4.6: ~{\sl Let $f_M$ be defined as above
for $M\in\cat$.  Then both $Im (f_M)\subset\prod_1^n\spy (M)$ and the
complement
$$\prod^n\spy (M)-(Im (f_M))$$ 
are $\party$ spaces of the first kind.}

\noindent{\sc Proof:} To verify ${\bf P}2$ for the case $\party
(M)=Im(f_M)$ notice that for $A\cap B=\emptyset$ in $M$, we have that
$\spy (A\sqcup B)=\spy (A)\times\spy (B)$ and hence $\spy (A\sqcup
B)^m=\spy (A)^n\times\spy (B)^m$. This then gives
\begin{eqnarray*}
   f_{A\sqcup B}(\spy (A\sqcup B)^m)&=& f_{A\sqcup B}(\spy (A)^m)\times
   f_{A\sqcup B}(\spy (B)^m) \\
   &=&f_A(\spy (A)^m)\times f_B(\spy (B)^m)
\end{eqnarray*}
and {\bf P}2 for this case follows. 

Notice at this point that 4.4 implies the existence of a {\sl relative} map
$$f_{M,N}:\prod^m\spy (M,N)\lrar\prod^n\spy (M,N).$$
and we can identify $\party (M,N)$ with $(Im f_{M,N})$ in the first case and
with its complement in $\spy (M,N)^{\times n}$ in the second. 

We now verify {\bf P}2 for spaces of the form $\party (M) =
\prod^n\spy (M)-(Im (f_M))$. Let $A,B\subset M$ as before, then
\begin{eqnarray*}
   &&\party (A)\times\party (B)\\
   &&= (\spy (A)^n-(Im f_A))\times (\spy (B)^n-(Im f_B))\\
   &&= \spy (A)^n\times\spy (B)^n- \left(\spy (A)^n\times Im f_B \cup
   Im f_A\times\spy (B)\cup Im f_A\times Im f_B\right)\\
   &&=\spy (A\sqcup B)^n - (Im f_A\times Im f_B) = 
   \spy (A\sqcup B)^n - (Im f_{A\sqcup B})\\
   &&=\party (A\sqcup B).
\end{eqnarray*}
The proposition follows.
\hfill\za

\noindent{\bf Corollary} 4.7: Let $f_M:\spy (M)^m\rightarrow\spy
(M)^n$ be as before, then the quotient monoids below form particle
spaces
$$\party (M) = \spy (M)^n//Im f_M~~\hbox{and}~~\party (M)=\spy (M)^n//
\left(\spy (M)^n-(Im f_M)\right).$$

\noindent{\sc Example} 4.8: Consider the diagonal map
$$M\fract{\Delta}{\ra 2}M\times M\fract{+}{\ra 2}\sp{2}(M).$$ We
extend it multiplicatively to a map $f_M:\spy (M)\rightarrow\spy (M)$
and it is direct to see that the complement of $(Im f_M)$ is $\config
(M)$.

\noindent{\sc Example} 4.9: Consider the map 
$$M\times M\lrar \spy (M)^{\times 3},~(a,b)\mapsto (a,b, a+b)$$
and extend it additively to a map $f_M:\spy (M)^{\times 2}\rightarrow
\spy (M)^{\times 3}$. Then $Im f_M$ corresponds to 
triples of configurations $(\zeta_1,\zeta_2,
\zeta_3)$ such that $\zeta_3=\zeta_1+\zeta_2$.

\noindent{\sc Example} 4.10: Spaces of pairwise disjoint configurations;
$DDiv^n$ (already defined in the introduction) can be described along the
lines formulated above. Assume for example $n=3$, then $DDiv^3(M)$ is the
complement in $\spy (M)^3$ of $(Im f_M)$ where $f_M$ is given by
$$
f_M:\spy (M)^3\lrar\spy (M)^3,~(\zeta,\eta,\psi)\mapsto
(\zeta+\eta, \zeta+\psi, \eta+\psi ).
$$


\vskip 10pt
\noindent{\Large\S5 Some Topological Properties}

Naturally $\party (M)$ inherits its topology from $\prod\spy (M)$
and the topology on $\spy (X)$ is the weak topology relative to the
subspaces $\sp{r}(X), r\geq 1$; that is a set $U\subset\spy (X)$ is
closed if and only if $U\cap\sp{r}(X)$ is closed for all $k$.

\noindent{\bf Lemma} 5.1:~{\sl Let $\party$ be a particle functor and
let $M$ be a manifold of dimension $n\geq 1$ such that $\party
(M)\neq\emptyset$. Then $\party (A)\neq\emptyset$ of all open
$A\subset M$.}

\noindent{\sc Proof:} Let $(\zeta_1,\ldots, \zeta_k)\in\party (M)$.
Then $S=\{\zeta_1\cup\cdots\cup\zeta_k\}$ is a finite set of points
and so there is always an injection of $\tau: S\rightarrow A$.  Since
$\party$ is a functor from $\cat$ to $\cat$, it follows that there is
an induced injection sending $\{(\zeta_1,\ldots,
\zeta_k)\}\in\party (S)$ into $\party (A)$ and the lemma follows.
\hfill\za

Recall that $\party$ is a self-functor of the category $\cat$ of
spaces and injections as morphisms. In particular, $\party$ takes
inclusions to inclusions.  Using the isotopy properties of $\party
(-)$ the following is not hard to establish.

\noindent{\bf Lemma} 5.2:~{\sl Let $M$ be compact with boundary and
denote by $M^{int}$ its interior. Then we have a homeomorphism
$\party (M)\cong\party (M^{int})$.}

\noindent{\sc Definition} 5.3: We let $\cat_n\subset \cat$
consist of the subcategory of $n$ dimensional ($n\geq 1$), smooth,
connected and compact manifolds.

\noindent{\bf Isotopy and injective homotopy}~: From the functorial
properties of $\party$, it is clear that any injective homotopy $h_t:
U\lrar M$; i.e. a homotopy through injective maps, induces a homotopy
of particle spaces; $\party (h_t):\party (U)\lrar\party (M)$.  An
isotopy between $f,g:N\rightarrow M$ is on the other hand a
differentiable homotopy through embeddings. It is ambiant if there is
an isotopy $F:M\times I\rightarrow M$ such that $F(x,0)=f(x)$ and
$F(x,1)=g(x)$.

Given any two points $p$ and $q\in int (M)$, $M$ connected, any smooth path
between them gives rise to an isotopy from $p$ to $q$. This isotopy
can be extended to an ambiant isotopy [Ko]. And generally one has

\noindent{\bf Lemma} 5.4:~{\sl On a manifold $M\in\cat_n$, 
there is an ambiant isotopy taking any finite set of interior points 
to any other set of interior points with the same cardinality.}

\noindent{\sc Proof:}
Two isotopic embeddings, via an isotopy $F:N\times I\rightarrow M$, 
need not be ambiant isotopic.
However there is an extension theorem of Thom that gives sufficient conditions
for when this is possible; namely when $N$ is compact and $M$ is closed.
A close inspection of the proof shows that the ambiant isotopy can
be chosen so that it leaves all points outside a compact neighborhood $V$ of
$N\times I$ fixed; [Mi]. 

Therefore and as long as this neighborhood $V$ misses
the boundary of $M$, the theorem of Thom still applies for non-closed $M$,
namely for $int(M)$. In our case, $N$ is a collection of points
$\{x_1,\ldots, x_m\}\in int (M)$ and hence is
compact. Let $\{y_1,\ldots, y_m\}\in int (M)$ be any other set
of $m$ points and choose paths $\gamma_i$ between $x_i$ and $y_i$ that lie
in the interior. Traveling along the paths (at different speeds if need
be in order to avoid intersections at any given time) gives an isotopy
$F:N\times I\rightarrow M$. By Thom's theorem, $F$ extends to
an ambiant isotopy and the lemma follows.
\hfill\za

\noindent{\bf Corollary} 5.5:~{\sl Let $N\in\cat_n$ be a connected
space and assume $\party (N)\subset\spy (N)^k$.  Then $\party (N)$ has
$\bbz^{+}\times\ldots\times\bbz^{+}$ components obtained as the
intersection of $\party (N)$ with
$\sp{n_1}(N)\times\cdots\times\sp{n_k}(N)$ for all tuples of positive
integers $(n_1,\ldots ,n_k)$. (Compare 3.5)}

\noindent{\sc Notation}:
In the case of particle spaces of the first kind, we can then index 
the components as follows
$$\party_{m_1,\ldots ,m_k}(M) = \party (M)\cap 
\sp{m_1}(M)\times\ldots\times\sp{m_k}(M)\subset\prod^k\spy (M),~m_i>0.$$
More generally, if $\party (M)$ is any particle space
given as a quotient $q:\party_1(M)\rightarrow\party (M)$ for some first
kind $\party_1$, then we define
$$\party_{m_1,\ldots ,m_k}(M)=q(\party_{1_{m_1,\ldots, m_k}}(M)).$$
We will see later (9.17) that the multidegrees $(m_1,\ldots, m_k)$ 
parametrize maps from $H_n(M;\bbz )$ into $H_n(\party (S^n,*);\bbz)$.

\noindent{\bf Lemma} 5.6:~{\sl Let $M\in\cat_n$, 
and let $N\subset M$ be an absolute neighborhood retract.
Then $\party (M,N)$ is connected.}

\noindent{\sc Proof:}
$N$ being as above, there is an open $U\subset M$ containing $N$ and
retracting to it via a retraction $r$. We assume this retraction is
injective on $N-U$ (think of a collar).  Given a multiconfiguration
$\{\zeta_1\cup\cdots\cup\zeta_k\}$ in $\party (M,N)$ (see note preceding
3.8), we let its support be the set of points making up the $\zeta_i$'s.
If this support lies in $U$, then the retraction $r$ takes 
$\{\zeta_1\cup\cdots\cup\zeta_k\}$ to $N$ and hence to basepoint in
$\party (M,N)$. Generally if ${\vec\zeta} = \{\zeta_1\cup\cdots\cup\zeta_k\}$
has support in $M-N$, then there always is an isotopy taking $\vec\zeta$
to an element $\vec\zeta\prime$ in $U$ (by lemma 5.4). Composing this with
$r$ gives at the end a path connecting $\{\zeta_1\cup\cdots\cup\zeta_k\}$
to basepoint and the lemma follows.
\hfill\za

\noindent{\bf Example} 5.7: Choose a basepoint $*\in M\in\cat_n$ which is
an interior point. Then $\party (M,*)$ is connected. We show in \S9
that if $M$ is $n$ connected then so is $\party (M,*)$.


\vskip 10pt
\noindent{\bf\Large\S6 Particle Spaces and Cofibrations}

\noindent{\bf\S6.1 Restrictions and Relative Constructions}:~
Fix a particle functor $\party$ and let $M\in\cat$ and $*\in N\subset
M$ closed.  Naturally $\spy (N)$ is a submonoid of $\spy (M)$ and
we define $\spy (M,N)$ as the quotient monoid $\spy (M)/\spy (N)$.
When $N=*$, it can be checked that this construction agrees with the
previously defined $\spy (M,*)$ in 2.5.
Now suppose $\party (-)$ is a functor of the first kind, then we
define
$$\party (M,N)=\left\{\zeta\in\spy (M,N)^k~|~ \zeta\cap (M-N)\in\party
(M-N)\right\}.
$$ 
If $\party (-)$ is obtained as the quotient of $\party_1(-)$ for some
particle functor of the first kind, then $\party (M,N)$ is 
obtained as a pushout construction
$$\matrix{\party_1 (M)&\lrar&\party_1(M,N)\cr
\downarrow&&\downarrow\cr
\party (M)&\lrar&\party (M,N).\cr}$$
In words, $\zeta\in\party (M,N)$ if 
$\zeta\cap (M-N)\in\party (M-N)$ with the additional contraint that as points
of $\zeta$ tend to $N$ they get identified with basepoint.

\noindent{\bf Remark} 6.1: Notice that $\party (M,N)$
has a canonical basepoint $\zeta = \langle *,*,\ldots \rangle$.
Observe as well that $\party (M,N)\simeq \party (M/N,*)$ (compare 2.6).

\noindent{\bf Lemma} 6.2:~{\sl Let $M\in\cat_n$, $N\subset M$. 
Then we have a quotient map
$\pi: \party (M)\lrar\party (M,N).$
If $N$ has boundary $\partial N$, we get a {\sf restriction}
$$r: \party (M)\lrar\party (N,\partial N).\leqno{6.3}$$}

\noindent{\sc Proof:} We simply need mention that $r$ is a special
case of $\pi$ as applied to the quotient $M\rightarrow M/(M-N)$ and
one can check that $\party (M, M-N)=\party (N,\partial N)$.
\hfill\za

\noindent{\bf Remark} 6.4: We can give an explicit description of
$\pi$ as follows. Let ${\vec\zeta}\in\party (M)$. Then since $\party
(M) = \party (N)+\party(M-N)$, we can write ${\vec\zeta} =
{\vec\zeta}_N + {\vec\zeta}_{M-N}$ where ${\vec\zeta}_N\in\party (N)$
and ${\vec\zeta}_{M-N}\in\party (M-N)$.  The correspondence
$${\vec\zeta}\mapsto {\vec\zeta_{M-N}}$$
is not continuous. However when post-composed with the quotient map
$$
\party (M-N)\lrar\party (\overline{(M-N)},\partial\overline{(M-N)})
\cong\party (M,{\overline N})=\party (M,N)
$$
it becomes so, hence yielding 6.2 (here we use the fact that $\party
(M-N)$ is homeomorphic to $\party (\overline{M-N})$).  On the
other hand, the correspondence ${\vec\zeta} = {\vec\zeta}_N$ yields
the {\sl restriction} map $r: \party (M)\lrar\party (N,\partial N).$

\noindent{\bf Remark} 6.5: There are different other restriction maps.
For instance, let $M_0\subset N\subset M\in\cat_n$,
then we have maps as follows
$$\party (M,M_0)\lrar\party (N,\partial N\cup M_0).$$
Notice also that given any morphism of pairs in $\cat$,
$(M,N)\hookrightarrow (M',f(N))$ we get an induced morphism
$$\party (M ,N)\lrar\party (M',f(N)).$$

\noindent{\bf\S6.2 Behaviour with respect to cofibrations}~:
From now on we restrict attention to the subcategory $\cat_n$, and hence
$\party:\cat_n\rightarrow\cat$. Associated to any pair $(M,N)\in\cat_n$, 
($N\subset M$ is of codimension $0$), we have the cofibration sequence
$$N\hookrightarrow M\rightarrow M/N.$$
Using the covariance of $\party$ with respect to
inclusions and using the restriction map constructed early in this
section, we can
apply $\party$ to the above sequence and get
$$\party (N)\lrar \party (M)\lrar \party (M,N).$$
More generally, we can start with the cofibration sequence
$$(N, N\cap M_0)\lrar (M,M_0)\lrar (M, N\cup M_0).$$
Then the following is a generalization of proposition 2.2, and
[B] (p:178);

\noindent{\bf Proposition} 6.6:~{\sl Consider the cofibration sequence
$(N, N\cap M_0)\lrar (M,M_0)\lrar (M, N\cup M_0)$ with $N\subset
M\in\cat_n$, $M_0\subset M$. Suppose that $M_0\cap N\neq\emptyset$, then
$$\party (N, N\cap M_0)\lrar \party (M,M_0)\fract{\pi}{\lrar}
\party (M, N\cup M_0)$$
is a quasifibration.}

\noindent{\sc Proof:} 
The submanifold $N\subset M$ being proper and compact, it has non empty
boundary $\partial N$ which we can assume wlog to be connected.
$\partial N$ has a tubular neighborhood $U_{\partial}\subset M$
when restricted to either $N$ or $M-N$ looks like a collar. Let
$U=N\cup U_{\partial}$, then there is
an {\sf isotopy} retraction of $r_t:U\lrar N$ which leaves $M-U$ and $N$ 
invariant ([Ko], chap.3). 
Consider at this point the subspaces 
$$
X_{k_1,\ldots, k_n}:=\left\{{\vec\eta}\in
\party (M, N\cup M_0)~|~{\vec\eta}\cap (M-N\cup M_0)\in
\party_{i_1,\ldots, i_n}(M - N\cup M_0),~i_j\leq k_j\right\}$$
(here $n$ is determined by $Par$)
and consider the open sets in $\party (M, N\cup M_0)$
$$
U_{k_1,\ldots, k_n}=\{(\zeta_1,\ldots, \zeta_r)\in X_{\vec k}~|~
(\zeta_1,\ldots, \zeta_r)~\hbox{contains
a non-empty subtuple in}~\party (U)\}.
$$
Write $X_{\vec k} = X_{k_1,\ldots, k_n}$ and similary $U_{\vec k}= U_{k_1,
\ldots ,k_n}$.
By construction, we have 
the following inclusions 
$$X_{<{\vec k}}: = \bigcup_i X_{k_1,\ldots, k_i-1,\ldots, k_n}
\subset U_{\vec k}\subset X_{\vec k}$$
and it is easy to see that over $X_{\vec k}-X_{<\vec k}$
the map $\pi:\party (M,M_0)\rightarrow\party (M,M_0\cup N)$ is
a direct product $(X_{\vec k}-X_{<\vec k})\times \party (N,N\cap M_0)$.
This direct product structure is actually given as follows: To ${\vec\zeta}
\in X_{\vec k}-X_{<\vec k}$ and ${\vec\eta}\in\party (N,N\cap M_0)$ we
associate the multiconfiguration ${\vec\zeta}+{\vec\eta}\in\pi^{-1}
(X_{\vec k}-X_{<\vec k})\subset\party (M, M_0)$. The sum 
${\vec\zeta}+{\vec\eta}$ is well defined since the support of $\vec\zeta$
is in $M-N\cup M_0$ and the support of $\vec\eta$ is in $N$.
This yields the first half of criterion 2.3. It is left to analyze
the attaching maps associated to this filtration.

One now sees that the isotopy retraction $r_t$ moves $\partial N$ away from
itself; $r_1(N)\subset N$, and squeezes the collar $U_{\partial}$ into $N$. 
This is done through a homotopy that is injective on $M-U_{\partial}$
and so from earlier considerations 
it induces a retraction at the level of particle spaces 
$r:U_{\vec k}\lrar X_{<{\vec k}}.$
Let ${\vec x} = (\eta_1, \ldots, \eta_r)\in U_{\vec k}$ and write
$$
(\eta_1+\zeta_1,\ldots,
\eta_r+\zeta_r)\in\pi^{-1}({\vec x}),
$$
where the $\zeta_i$ are in $\party (N, N\cap M_0)$.
Since $(\eta_1, \ldots, \eta_r)\in U_{\vec k}$, there exist a
subtuple $(D_1,\ldots, D_r)\subset (\eta_1, \ldots, \eta_r)$ such that
$D_i\in U-N$. The retraction $r$ moves the $D_i$ inside $N$ (so at time
$t=1$, $r_1(U)\subset N$) and hence at
the level of preimages we have a lifting
\begin{eqnarray*}
\party (N, N\cap M_0)&\fract{{\tilde r}_1}{\ra 2}&\party (N, N\cap M_0)\\
(\zeta_1,\ldots, \zeta_r)&\mapsto&
\left(r_1(\eta_1)+r_1(D_1),\ldots, r_1(\eta_r)+r_1(D_r)\right)
\end{eqnarray*}
where ${\tilde r}_1$ is the attaching map 
$${\tilde r}_1:\party (N,N\cap M_0)\fract{\cong}{\lrar}
\pi^{-1}({\vec x}){\ra 2} \pi^{-1}(r_1({\vec x}))\fract{\cong}{\lrar}
\party (N,N\cap M_0)\leqno{6.7}$$
Since $M_0\cap N\neq\emptyset$, we use lemma 6.2 to deform
$(r_1(D_1),\ldots, r_1(D_r))$ through a path in $\party (N)$ to
a tuple in $M_0\cap N$ and hence to basepoint in $\party
(N, M_0\cap N)$. This produces a homotopy inverse for
${\tilde r}_1$ and the proposition follows by 2.2.  \hfill\za

\noindent{\bf Remark} 6.8: 
When $N\cap M_0$ is empty, then $\party (N)$ does generally split into
components. In this case, the attaching map ${\tilde r}_1$ 
switches components and it has no homotopy inverse. We deal with this
case in the next section.

\noindent{\bf Remark} 6.9: The proposition is also not true if $N$ is
not of codimension $0$ in $M$. For example, let $\party = \config$ be
the functor of disjoint unordered points (see introduction). We show
in 11.1 that $\config (D^n,\partial D^n)\simeq S^n$. Suppose in this
case that $M=D^n$, $N=\{(x_1,\ldots, x_{n-1},0)\}\subset D^n$ is an $n$-th
face and let $M_0=\partial D^n-N$. If 6.6 were to apply in this case, 
then we get a quasifibering
$$\config (D^{n-1},\partial D^{n-1})\lrar\config (D^n,M_0)\lrar
\config (D^n,\partial D^n).$$
But since $\config (D^n,M_0)$ is contractible, we would have proved that
$S^{n-1}$ is weakly homotopy equivalent to $\Omega S^n$ which is obviously
false.

\noindent{\bf\S6.3 Connectivity properties}~:
At this point, we would like to analyze the connectivity of the spaces
$\party (M,*)$. The following is a direct consequence
of 6.6 above.

\noindent{\bf Proposition} 6.10:~{\sl Let $M$ be any $n-1$ connected
finite CW complex ($n>1$). Then $\party (M,*)$ is also $n-1$ connected.}

\noindent{\sc Proof:} The proof is a standard induction on cells of
$M$. First since $M$ is $n-1$ connected, it has a CW decomposition
with cells starting in dimension $n$ attaching to a basepoint $*$.
Let $M^{(i)}$ denote the $i$-th skeleton of $M$ (here of course $i\geq
n$).  The inclusion of $M^{(i)}$ into the next skeleton gives a
cofibration sequence
$$M^{(i)}\lrar M^{(i+1)}\lrar\bigvee S^{i+1}$$
which yields by 6.6 a quasifibration 
$$\party (M^{(i)},*)\lrar\party (M^{(i+1)},*)\lrar\prod
\party (S^{i+1},*),~i\geq n.\leqno{6.11}$$
Suppose that $\party (S^n,*)$ is $n-1$ connected, then $\party
(M^{(n)},*)=\prod\party (S^n,*)$ is also $n-1$ connected and the the
long exact sequence in homotopy attached to 6.11 shows that $\party
(M^{(n+1)},*)$ is $n-1$ connected as well. Proceeding inductively,
we can establish the claim as soon as we show that
$\party (S^n,*)$ is $n-1$ connected. This is done in 9.2.
\hfill\za


\vskip 10pt
\noindent{\bf\Large\S7 ``Partial'' Group Completion}

Start with the simplest particle space $\spy (M)=\coprod\sp{n}(M)$ and
notice that
$$\spy (M)^+\simeq\bbz\times\spy (M,*)$$
where $\spy (M)^+$ means group completion with respect to the monoid
structure. One way of constructing $\spy (M)^+$ is to consider 
the space of {\em infinite} configurations of the form $\sum n_ix_i$
where the $n_i$ are not necessarily positive.
Choose a sequence $\{z_i\}_{i=1}^{\infty}$ of points in $M$, and
write $\eta = \sum z_i$ for the corresponding infinite configuration.
Then $SP(M)= \spy (M)^+$ is
the set of infinite configurations $\zeta$ such that
$\zeta-\eta$ is a finite but not necessarily positive configuration.

It is our desire to construct an analogue of $SP(-)$ (which we denote
by $Par(-)$) for the more general particle functors $\party (-)$.  The
end result would be some sort of ``group completion'' with respect to
the partial monoid structure of $\party (M)$. This is done in a very
standard way.

As always, let $M$ be compact (connected) and $A\subset M$ a closed
non-empty ANR (typically $A = \partial M$ for example).  We can
``stabilize" $\party (M)$ by ``marching toward A". Let $U$ be a
tubular neighborhood of $A$ which we assume to retract to $A$ via a
retraction $r$ which is injective outside of $U$.  Let $U_{i\in\bbz^+}$ be a
nested sequence $U_{i+1}\subset U_i\subset U$ The $U_i-U_{i+1}$ being
open, $\party (U_i-U_{i+1})\neq\emptyset$ according to lemma 5.1, and
so we can choose $\vec\eta_i\in \party (U_i-U_{i+1})$. We choose
$\vec\eta_i$ to be ``minimal'' in the sense that no smaller subtuple
of it lies in $\party (U_i-U_{i+1})$.  Now notice that we have an
inclusion given by summing with $\vec\zeta_i$ in the partial monoid
structure on $\party (M-U_{i+1})$;
$$\party (M-U_{i})\fract{+\vec\zeta_i}{\ra 3}\party (M-U_{i+1})$$
We can now make the definition

\noindent{\sc Definition} 7.1: For $M$, $A$ and $U$ as above, we define
$$\displaystyle Par (M) = \lim_{\vec\zeta_i}\left(\party (M-U_i)
\fract{+\vec\zeta_i}{\lrar}\party (M-U_{i+1})\right).$$

\noindent{\bf Remark} 7.2: It should be clear that $Par(M)$ doesn't
depend (up to homeomorphism) on the choice (up to isotopy) of the
stabilizing sequence $\vec\eta_i$ or of the nested sequence $\{U_i\}$.
It is equally clear (for the same reasons given earlier for the case
$\party=\spy$) that components of $Par(M)$ are homeomorphic.

\noindent{\bf Remark} 7.3: We can define $Par(M,N)$ for pairs $(M,N)$
by taking suitable direct limits over $\party (M-U_i,N-U_i\cap
N)$. When $A = \partial M$ for example (or a subset of it),
$N\cap\partial M\neq\emptyset$, then we can stabilize with respect to
a sequence of multiconfigurations $\{\eta_i\}$ converging to a point
$p\in N\cap\partial M$. By a homotopy (again injective in the complement
of $U$) we can retract points of $\eta_j$ to $p$ and this shows (in this
case) that
$$Par(M,N)\simeq\party (M,N).$$

\noindent{\bf Remark} 7.4: One may observe that if $p\in\partial
M\neq\emptyset$,then we can stabilize with respect to a tuple
${\vec\eta}_i$ converging to $p$ (that is the
sequence of points making up each $\eta_i$ converges to $p$). In this
case one can show that
$$Par_c(M)\simeq\party (M,p)$$

\noindent{\sc Example} 7.5: We consider the case $\config
(M)\subset\spy (M)$. We assume $M$ has boundary $\partial
M\neq\emptyset$ and stabilize as above with respect to a collar $U$ of
$M$ by choosing a sequence of {\sl distinct} points $z_i$ marching
towards $\partial M$. Addition of points $z_i$ yields maps and commutative
diagrams
$$\matrix{
   C_i(M)&\fract{+z_i}{\lrar}&C_{i+1}(M)\cr
   \decdnar{\subset}&&\decdnar{\subset}\cr
   \sp{i}(M)&\fract{+}{\lrar}&\sp{i+1}(M).
\cr}$$ 
which in the limit yield a map $C(M)\lrar SP(M)$. Now it isn't hard
to see that $C(M)$ breaks into $\bbz$ components which we write
$\config (M,p)$ for some $p\in\partial M$ (see 7.4). We then get a map
$$\alpha: \config (M,p)\lrar\spy (M,p)$$
It is interesting to notice that in the case $M=\bbr^n$ for example,
the left hand side is equivalent in homology to a component of
$\Omega^nS^n$ (cf. \S12) and hence $\alpha_n$ is homologous
to the map
$$\Omega_0^nS^n\lrar\spy (\bbr^n,*)\simeq *$$
obtained by looping $n$ times the natural inclusion $S^n\lrar\spy (S^n,*)$.

\noindent{\sc Definition}: A morphism $\iota: Par (N)\lrar Par (M)$
will mean an inclusion 
such that $\iota (\eta_N) =\eta_M$. The
following is a generalization of 6.6

\noindent{\bf Theorem} 7.6:~{\sl Given a cofibration sequence $N\lrar
M\lrar (M, N)$, $N,M\in\cat_n$, and assume $\partial
N\neq\emptyset$. Then
$$Par (N)\fract{i}{\lrar} Par (M)\fract{\pi}{\lrar}\party (M, N)$$ 
is a homology fibration.}

\noindent{\sc Proof:} First we convince ourselves that the preimage
under $\pi$ is homeomophic to $Par(N)$ and so $i$ is contructed
up to homeomorphism. The proof now amounts to
showing that the attaching map ${\tilde r}_1:Par (N) \rightarrow Par
(N)$ is a homology equivalence. Recall (\S6) that ${\tilde r}_1$ is
obtained by moving particles of $N$ away from a collar $U$ of
$\partial N$ and then adding a given element ${\vec\nu}\in \party
(r_1(U-N))$.

Let ${\vec\eta}_i$ be a stabilizing sequence used in constructing
$Par(N)$ with respect to some tubular neighborhood of $\partial N$. 
Then up to isotopy, $\vec\eta$ can be
written as a finite sum over some ${\vec\eta}_j$'s. To see this, we 
observe that via an isotopy (if necessary) we can bring points of
${\vec\eta}_1$ (for example) to points of $\vec\eta$ and hence
${\vec\eta}-{\vec\eta}_1$ will be positive (by minimality of $\eta_1$)
and belongs to $\party (N)$.  Reiterating this argument shows that
${\vec\eta} = \sum{\vec\eta}_{j}$ (a finite sum). By construction of
$Par(N)$ as a direct limit over addition of the
$\vec\eta_j$'s, $+\vec\nu$ necessarily induces a homology isomorphism
and the claim follows.  \hfill\za

\noindent{\bf Remark} 7.7: We emphasize again that $N$ needs to be of
codimension 0 in $M$ (see remark 6.9).  Here's (another) example where
7.6 doesn't hold if $N$ not of zero codimension.  Consider the
particle space $Div^2(M)\subset\spy (M)\times\spy (M)$ consisting of
pairs of configurations on $M$ with no points in common and let
$N=*\hookrightarrow M$. Then $Div^2(N)=Div^2(*)=\emptyset$ and there
isn't much sense to the sequence $Div^2(*)\rightarrow
Div^2(M)\rightarrow Div^2(M,*).$ However if we replace $*$ by the open
disc $D^n$, then $Div^2(D^n)$ has now reasonable properties; in fact
it is given by
$$\Omega^n(K(\bbz,n)\vee K(\bbz,n))~~(cf.\S12)$$
and $Div^2(D^n)\lrar Div^2(M)\rightarrow Div^2(M,*)$ is a quasifibration.


\vskip 10pt
\noindent{\bf\Large\S8 Scanning Smooth Manifolds}

The term ``scanning" is borrowed from Segal ([S1]).  We assume
as usual that $M$ is smooth and connected $n$-dimensional manifold.

\noindent{\bf\S8.1. Scanning parallelizable manifolds:~}
The scanning process is best pictured when $M$ is parallelizable
(i.e. $M$ has trivial tangent bundle). Examples of such manifolds are
Lie groups or any oriented three dimensional manifold.  Without loss
of generality, we restrict attention below to $\party (M) = \spy (M)$
(the more general situation is treated in the exact same way).

\noindent{\bf Definition} 8.1: Let $M^n$ be as above and let $N$ be a
closed subset of $N$. The {\em pair} $(M,N)$ is said to be {\em
parallelizable} if $M-N$ is parallelizable.  $M$ is stably
parallelizable if $(M,*)$ is parallelizable for instance.  Riemann
surfaces are examples of stably parallelizable surfaces, as well as
compact, oriented, spin four manifolds.

Put a metric on $M$ and consider the unit disc bundle $\tau M$ lying
over $M$. Let's assume for now that $\partial M=\emptyset$.)
Via the exponential map we can identify a neighborhood of
every point $x\in M$ with the fiber at $x$.  Denote such a
neighborhood by $D(x)\subset M$. When $M$ is parallelizable, the
fibers over $\tau M$ are canonically identified with a disc $D^n$ and
hence one can identify canonically the pairs $(\bar D(x),\partial
\bar D(x))$ for every $x\in M$ with $(D^n,\partial D^n) =
(S^n,\infty)$ (where the north pole $\infty$ is chosen to be the
basepoint in $S^n$.)

Given a configuration $\zeta\in\sp{d}(M)$ and an $x\in M$, then
$\zeta\cap D(x)$ is a configuration on $D(x)$ and its image 
under the restriction map 
$$\spy (D(x))\lrar \spy (\bar D(x),\partial 
\bar D(x))= \spy (S^n,\infty )$$ 
is denoted by $\zeta_x$. Notice that the correspondence $\zeta\mapsto
\zeta_x$ is now continuous (while the correspondence $\zeta\mapsto
\zeta\cap D^n(x)$ was not to begin with). Starting with $\zeta\in\sp{d}(M)$,
we hence get a map
$$
S_d:\sp{d}(M)\lrar\map{d}(M,\spy (S^n,*)),~~\zeta\mapsto f_{\zeta}:
f_{\zeta}(x) = \zeta_x.
$$
The scanning map $S$ is now given as $\sqcup S_r$.

\noindent{\bf Scanning manifolds with boundary}:~ We're in the case
$\partial M\neq\emptyset$. We can still scan the interior $M-\partial
M$ and alter the topology as points tend to $\partial M$.

Consider the open interior $M^{int}=M-\partial M$ and
let $\sp{r}_{\gre}(M^{int})$ be the subspace of
$\sp{r}(M)$ consisting of configurations of points that are at least
$2\gre$ away from the boundary $\partial (M)$. Choose
$\zeta\in\sp{r}_{\gre}(M^{int})$.  Then by scanning the
interior using discs of radius $\gre$, it is clear that $S_{\zeta}$
maps $x$ to basepoint for $x$ sufficiently near the boundary. This
gives rise to a map
$$ 
\sp{r}_{\gre}(M^{int})\lrar\map{r}(M/\partial M,\spy (S^n,*)).
$$ 
As $\gre\rightarrow 0$, one obtains in the limit a map
$$ 
S_r: \sp{r} (M)\lrar\map{r}(M/\partial M,\spy (S^n,*))
\simeq\map{r}(M/\partial M, K(\bbz, n))).  
$$ 
In exactly the same way, one obtains for each parallelizable pair
$(M,N)$ a map
$$ S: \spy (M-N)\lrar\map{}(M, N\cup\partial M,\spy (S^n,*))$$
where the right hand side consists of all based maps sending $N$
into basepoint $*\in\spy (S^n,*)$.

\noindent{\bf Scanning smooth manifolds}:~The general case of $M$ not
necessarilt parallelizable and $\party$ any particle functor is
treated similarly. The starting point is again the unit disc bundle
$\tau M$.  Compactifying each fiber yields a bundle
$$S^n\lrar\widehat{\tau M}\lrar M$$
to which we associate the bundle of configurations
$$\party (S^n,*)\lrar E_{\party}\lrar M.\leqno{8.2}$$
by applying $\party$ fiberwise.
Note that $\widehat{\tau M}$ has a ``zero" section $\nu$ sending each
$x\in M$ to the compactifying point in the fiber. We label this point
by $*$. Clearly, such a section extends to a zero section of $\party
(S^n,*)\lrar E\lrar M$ also denoted by $\nu$. We denote by
$\sect{} (M,A,\party (S^n,*))$ the space of sections restricting 
to $\nu$ on $A\subset M$.

The exponential map again provides a cover of $M$ by neighborhoods
$\bigcup_{x\in M} D^n(x)$ with respect to which we can scan. Cutting a
neighborhood $\bar D^n\subset M$ yields a cofibration $M-\bar D^n
\hookrightarrow M\rightarrow (\bar D^n,\partial \bar D^n)$ and hence
we get``retriction" maps
$$\pi_x: \party (M)\lrar \party (\bar D^n(x),\partial \bar
D^n(x)),~\forall x\in M$$ 
Starting with an element in $\party M$, one can restrict via $\pi_x$
to neighborhoods as in 8.1.  The elements of $\map{}(D^n(x),\party
(S^n,*)$ are now local sections of 8.2 and one gets the correspondence

\noindent{\bf Lemma} 8.3: Let $M\in\cat_n$ and $N\subset M$ a closed ANR. 
Then scanning yields a map
$$\party (M-N)\lrar\sect{}(M,N\cup\partial M,\party (S^n,*)).$$

\noindent{\bf\S8.2 The Electric Field Map}: The scanning map is closely
associated to another map of Segal ([S1]).
Given $k$ points in $\bbr^n$, we can electrically charge them and
hence we can associate to them an electric field $E$ which is a
function on $\bbr^n$ taking values in $\bbr\cup\infty$ (where $\infty$
is reached at the charge points). Since the electric field intensity
decays away from the charges, we can extend the previous function to
$\bbr^n\cup\infty$ by taking $\infty$ to the basepoint $0\in S^n$. So
to $k$-points $\{q_1,\ldots, q_k\}$ in $\bbr^n$, we associated the map
$E_{q_1,\ldots, q_k}:S^n\rightarrow S^n$ which is based and of degree
$k$ (since$E^{-1}(\infty)=\{q_1,\ldots, q_k\}$). This then shows that
$E_{q_1,\ldots ,q_r}\in \Omega_k^n(S^n)$ as desired. It is easy to see
that $S$ corresponds to $E_k$; that is

\noindent{\bf Lemma} 8.4:~{\sl The maps $S, E_k:
C_k(\bbr^n)\rightarrow\Omega_k^n(S^n)$ are homotopic.}

\noindent{\sc Proof:} Let $C_{\epsilon}(\bbr^n,k)$ be the subset of
divisors $D\in C(\bbr^n,k)$ such that $D=p_1+\cdots +p_k$, $p_i\neq
p_j$ and $|p_i-p_j|>\epsilon$. We consider the electric field map
$E_D$ associated to $D$. Let $B_{\epsilon}(x)$ be the ball of radius
$\epsilon$ around the point $x\in\bbr^n$. Then by shrinking the vector
field one can confine it inside the balls $B_{\epsilon}(p_i)$ so that
it is zero outside of these balls. If we choose the electric fields to
die out linearly, then the map $E_D$ is seen to correspond to scanning
the configuration $D$ with a ball of radius $\epsilon$. Now since
$C(\bbr^n, k)=\bigcup_{\epsilon>0}C_{\epsilon}$ the lemma follows.
\hfill\za


\vskip 10pt
\noindent{\bf\Large\S9 Proof of Main Theorems 1.1 and 1.3}

\noindent{\sc Notation}: Recall that $\party (S^n,*)$ has a
``preferred'' identity $*$. For each pair of spaces $(M,N)$, we will
write $\map{}(M,N,\party (S^n,*))$ for the space of continuous maps
from $M$ into $\party (S^n,*)$ which send $N$ to *. 

\noindent{\bf\S9.1 The Homology Equivalence 1.1}: A good starting
point is a quick review of the handle decomposition of a manifold
(\`a la Thom-Smale-Milnor): Let $M$ be a
given smooth, compact manifold of dimension $n$. Then there exists a
sequence of manifolds with boundary
$$M_0=D^n\subset M_1\subset\cdots\subset M_{n-1}\subset M$$
such that $M_i$ for $0<i<n$ is obtained from $M_{i-1}$ by attaching a
number of $i$-handles to its boundary. The manifold $M$ is now
obtained from $M_{n-1}$ by attaching an $n$-dimensional cell. If $M$
has boundary, then one sets $M_{n-1}=M$.

A handle $H^i$ of index $i$ is a copy of $D^n = D^i\times D^{n-i}$
which is attached to $\partial M_{i-1}$ via an embedding $h:
S^{i-1}\times D^{n-i}\hookrightarrow\partial M_{i-1}$. The resulting
space is a smooth manifold. The pair $(D^i,S^{i-1})$ is called the
{\em core} of the handle, while $S^{i-1}\times D^{n-1}$ is referred to
as the {\em attaching region}.

\noindent{\bf Proposition} 9.1:~~{\sl Let $D^n$ be the closed unit
disc. For $0<k \leq n$ we have homotopy
equivalences
$$\party (D^n, S^{k-1}\times D^{n-k})\simeq
\Omega^{n-k}\party (S^n,*),$$
while for $k=0$ we have a homology equivalence
$H_*(Par (D^n);\bbz )\simeq H_*(\Omega^n (\party (S^n,*));\bbz)$.}

\noindent{\sc Proof:} Let $H^k$ be a handle $D^n=D^k\times D^{n-k}$ of
index $k$ attached along $(\partial D^k)\times D^{n-k}$.  We write
$A^k = S^{k-1}\times D^{n-k}$ and the goal is then to show that
$Par(H^k,A^k)\simeq \Omega^{n-k}\party (S^n,*)$ for $0\leq k\leq n$
(note that $A^0=\emptyset$).

The proof uses the cofibration sequence described in [B]. Let
$I_k\subset D^n=[0,1]^n$ denote the subset of $(y^1,\ldots, y^n)$ such
that $y^i=0$ or $y^i=1$ for some $i=1,\ldots, k-1$, or $y^k=1$ (that
is $I_k$ consist of all the boundary faces of $D^k \subset
D^n=D^k\times D^{n-k}$ safe the face $y^k=0$). Now let
$H_k=[0,1]^{k-1}\times [0,{1\over 2}]\times [0,1]^{n-k}$.  Then there
is a cofibration sequence
$$
(H_k, H_k\cap I_k)\lrar (D^n, I_k)\lrar (D^n, H_k\cup I_k).
\leqno{9.2}
$$
The pair $(H_k, H_k\cap I_k)$ can be identified with $(D^n,
S^{k-2}\times D^{n-k+1})$ hence representing a $k-1$-handle
$(H^{k-1},A^{k-1})$, while the pair $(D^n, H_k\cup I_k)=(D^n,
S^{k-1}\times D^{n-k})$ represents a handle $(H^k,A^k)$ of index $k$.
Applying $\party$ to $9.2$ yields the quasifibration (theorem 6.6)
$$
\party (H^{k-1},A^{k-1})\lrar\party (D^n,I_k)\fract{\pi}{\lrar}\party 
(H^k,A^k).$$
The proof now proceeds by downward induction on $k$.  Observe first
that 9.1 is obviously true when $k=n$. Suppose it is true for some
$n\geq k>0$. Then one can consider the following diagram of
quasifibrations
$$\matrix{
   \party (H^{k-1},A^{k-1})&\lrar&\Omega^{n-k+1}\party (S^n,*)\cr
   \downarrow&&\downarrow\cr
   \party (D^n,I_k)&\lrar& PS \cr
   \decdnar{\pi}&&\decdnar{}\cr
   \party (H^k,A^k)&\lrar&\Omega^{n-k}(\party (S^n,*))\cr}
\leqno{9.3}$$
By induction the bottom map is a homotopy equivalence while the middle
map is also a homotopy equivalence since $\party (D^n,I_k)$ is
contractible whenever $k\geq 1$ (this follows from the fact that when
$k\geq 1$, $I_k\neq \emptyset$ and there is a retraction of $D^n$ onto
$I_k$ which is injective on the complement of a tubular neighborhood 
of $I_k$). It follows that the top
inclusion must be a homotopy equivalence whenever $k\geq 1$ (this establishes
the first claim in 9.1).

Remains to treat the case $k=1$. Since $A^0=\emptyset$ the left hand side
in 9.3 is not a quasifibration anymore and we have to pass to the $Par$ 
functor. We can then consider the diagram of fibrations
$$\matrix{
   F&\lrar&\Omega^{n}\party (S^n,*)\cr
   \downarrow&&\downarrow\cr
   Par (D^n,I_1)&\lrar& PS \cr
   \decdnar{\pi}&&\decdnar{}\cr
   \party (H^1,A^1)&\fract{\simeq}{\lrar}&\Omega^{n-1}(\party (S^n,*))\cr}
$$
where $F\simeq \Omega^{n}\party (S^n,*)$ is the homotopy fiber for the
l.h.s and $Par(D^n,I_1)\cong\party (D^n, I_1)$ (see 7.3) is contractible.
Since the l.h.s is a homology fibration (theorem 7.6),
the inclusion of the preimage $Par(D^n,A^0)=Par(D^n)$ into $F$ is a
homology equivalence and the proof is complete.  \hfill\za

\noindent{\bf Corollary} 9.4:~{\sl $\party (S^n,*)$ is $n-1$ connected.}

\noindent{\sc Proof}: We have that $\pi_k(\party
(S^n,*))=\pi_{0}(\Omega^{k}\party (S^n,*))$ and the latter is
trivial whenever $k<n$ since $\Omega^{k}\party (S^n,*)$ is identified with
$\party (D^{n-k}\times D^{k}, S^{n-k-1}\times D^{k})$ and the latter
is connected by lemma 5.6.
\hfill\za

\noindent{\sc Example}: When $\party = \config$ (again), 
we show in 11.1 that $\config
(S^n,*)\simeq S^n$ (this is an early observation of Segal and McDuff).
Naturally $S^n$ is $n-1$ connected and 9.1 shows that for $0\leq k < n$
$$\config (D^n, D^k\times S^{n-k-1})\simeq\Omega^k S^n.$$

\noindent{\bf Theorem} 9.5:~{\sl Let $M\in\cat_n$ and $N$ a closed ANR in
$M$. Suppose either $N$ or $\partial M$ non-empty. Then scanning
induces a homology equivalence
$$
S_*: H_*\left( Par(M-N); \bbz \right) \fract{\cong}{\ra 2}
H_*\left(\sect{}(M,N\cup\partial M, \party (S^n,*));\bbz\right).$$}

\noindent{\sc Proof:} Since $Par$ is an isotopy functor, then
$Par(M-N)=Par(M-T(N))$ where $T(N)$ is a tubular neighborhood of $N$
and so wlog we can assume that $N$ is of codimension 0. 
We consider the case $N\neq\emptyset$ and $\partial M=\emptyset$
(the other cases are treated similarly). Since $Par (M-N)=Par(M-int
(N))$, we can assume that $M-N$ is compact and has boundary $\partial
N$. So $M-N$ has a handle decomposition
$$
M_0=D^n\subset M_1\subset\cdots\subset M_{n-1}=M-N
$$
where the handles we attach have index at most $n-1$ and all the $M_i$
have boundary $\partial M_i\neq 0$. The proof now proceeds by
induction on $i$. Since the number of handles we attach at each stage
(finitely many) is immaterial for the arguments below we
might as well assume we're
only attaching one handle at a time. That is
$$
M_i = M_{i-1}\cup H^i,~\hbox{and}~\partial M_{i-1}\cap H^i=S^{i-1}\times 
D^{n-i}.
$$
Consider the following two cofibrations;
$$
M_{i-1}\lrar M_i\lrar (D^i\times D^{n-i},S^{i-1}\times D^{n-i}), i<n,
\leqno{9.6}
$$
and the one induced from the handle attachment
$$
(H^i, H^i\cap\partial M_i)\lrar (M_i,\partial M_i)\lrar (M_i,H^i\cup\partial
M_i)=(M_{i-1},\partial M_{i-1}).
\leqno{9.7}
$$
We now apply the
functor $\sect{}$ to $9.7$ and get the {\em fibration}
$$\small
  \sect{} (M_{i-1},\partial M_{i-1},\party (S^n,*))\rightarrow
  \sect{}(M_i,\partial M_i,\party (S^n,*)\rightarrow
  \sect{} (H^i, H^i\cap\partial M_i,\party (S^n,*)).
\leqno{9.8}$$
Since $E_{\party (H^i)}$ over $H^i$ is trivial, we can replace 
$\sect{} (H^i, H^i\cap\partial M_i,\party (S^n,*))$ by an iterated
loop space as follows
\begin{eqnarray*}
   \map{}(H^i, H^i\cap\partial M_i;\party (S^n,*))&=&
   \map{}(D^i\times D^{n-i},\partial D^{n-i};\party (S^n,*))\\
   &=&\map{}^*(S^{n-i}\times D^i; \party (S^n,*))\\
   &=&\Omega^{n-i}\party (S^n,*).
\end{eqnarray*}
On the other hand one can apply the functor $Par$ to 9.6 and obtain
a quasifibration which maps via scanning into 9.8 as follows
$$
\begin{array}{ccc}
   Par(M_{i-1})&\fract{S}{\lrar}&\sect{}(M_{i-1},\partial M_{i-1};
   \party (S^n,*))\\
   \downarrow&&\downarrow\\
   Par(M_{i})&\fract{S}{\lrar}&\sect{}(M_{i},\partial M_{i};\party (S^n,*))\cr
   \downarrow&&\downarrow\\
   \party(D^i\times D^{n-i},S^{i-1}\times D^{n-i})&\fract{\simeq}{\lrar}&
   \Omega^{n-i}\party (S^n,*).
\end{array}
\leqno{9.9}
$$
The bottom map is a homotopy equivalence whenever $1\leq i\leq
n$ by 9.1. When $i=1$, $M_{i-1}=D^n$ and the top map is a homology
equivalence (here $\sect{}(M_{0},\partial M_{0};
\party (S^n,*))$ is again identified with $\Omega^n\party (S^n,*)$.)
By a standard spectral sequence argument, it follows that the middle
map $Par(M_{1})\lrar\sect{}(M_{1},\partial M_{1};\party (S^n,*))$ is
a homology equivalence and the argument proceeds by induction.
\hfill\za

\noindent{\bf Remark} 9.10: The theorem above is not true if both $N$
and $\partial M$ are empty. In that case ($M$ is closed)  
$Par(M)$ is not even defined and it doesn't even
hold true that the components of $\map{}(M,\party (S^n,*))$ are
homotopy equivalent. For example, we show in [K2] that components of
$\map{}(M_g,\bbp^n)$ (which is a special case of 11.7) do differ in
homotopy type.

\noindent{\bf Theorem} 9.11:~
{\sl Let $N,M$ be as in 9.4 and suppose
$\pi_1(Par(\bbr^n))$ is abelian, then scanning is a homotopy equivalence
$$Par (M-N) \fract{\simeq}{\ra 2} \sect{}(M,N\cup\partial M, \party
(S^n,*)).$$}

\noindent{\sc Proof}: Consider 9.9 again and the case $i=1$.  When
$\pi_1(Par(\bbr^n))$ is abelian, the l.h.s in 9.9 becomes a
quasifibration (this is explained in 9.12 below) and hence the top map
is a weak homotopy equivalence. Since the spaces involved have the
homotopy type of CW complexe we get a homotopy equivalence and hence
an equivalence in the middle. The rest of the proof is obtained by
induction knowing that the bottom map is always a homotopy equivalence
when $1\leq i\leq n-1$.  \hfill\za

\noindent{\bf\S9.2 Good Functors}: 
The functor $Par$ is {\sl good} if it turns cofibrations $N\rightarrow
M\rightarrow M/N$, $N,M\in\cat_n$ into quasifibrations. A
straightforward examination of the proof of 9.5 shows that scanning
induces a weak homotopy equivalence
$$Par(M-N)\fract{\simeq}{\ra 2}\sect{}(M,N\cup\partial M, \party (S^n,*))$$
whenever $Par$ is good and $N$ and $M$ are as in 9.5
(note that the space of sections has the homotopy type of a CW complex;
cf. [BS], lemma 6.5). 
The condition needed in
theorem 9.11 is of course slightly weaker.

\noindent{\bf Lemma} 9.12:~{\sl The functor $Par$ is good if it 
abelianizes fundamental groups; that is if $\pi_1(Par (M))$ is abelian
for any $M\in\cat_n$.}

\noindent{\sc Proof:} We need show that $Par$
applied to cofibrations yields quasifibrations. This boils down to
showing that the attaching maps given by addition of particles (see
6.7) are homotopy equivalences.  These attaching maps which take the
form $Par(M)\fract{+\zeta}{\lrar} Par (M)$ are homology equivalences
for any twisted coefficients (by construction of $Par(M)$ as a direct
limit over these additions). When $\pi_1(M)$ is abelian, the map
$+\zeta$ induces an isomorphism of fundamental groups as well.  
This implies that the attaching maps must be a homotopy
equivalences and the lemma follows.  \hfill\za

\noindent{\bf Example} 9.13: Since $\pi_1(\spy (X))=H_1(X)$ (for any space
$X$) we automatically
have that $SP$ is a good functor. More is true in this case for one
can show that $\pi_1(\sp{n}(X))$ is already abelian when $n\ge 2$. To see
this, let $\alpha\in\pi_1(\sp{n} (X))$, and $q: X^n\rightarrow \sp{n}(X)$
for the quotient map. The loop
$\alpha$ (representing a class in $\pi_1$) can be homotoped away from
the branched points for $n>1$ (with basepoint * fixed) and hence it can be
lifted to $X^n$. Since $q^{-1}(*) = *$, it follows that $\alpha$ lifts
to a loop in $X^n$. The rest of the claim follows from this observation.

\noindent{\bf Example} 9.14: The functor $C$ is not ``good'' 
and so the homology equivalence of 9.5 for the case $Par=C$
cannot be upgraded in general to a homotopy equivalence. A standard
example is already provided by the closed unit disc $D^n$. Consider
$C(\bbr^n)\cong C(D^n)$ and hence from theorem 9.5 we infer that
$$H_*(C(\bbr^n);\bbz )\cong H_*(\Omega^nS^n;\bbz ).$$
At the level of components $C(D^n) =\bbz\times\config (D^n,*)$
where $\config (D^n,*)$ is as described in example 7.4. It is known that
$\pi_1(\config (D^n,*))\cong\Sigma_{\infty}$ for $n>2$ (and
is the braid group for $n=2$). Since $\pi_1(\Omega^nS^n )\cong\bbz_2$,
theorem 9.5 in this case couldn't possibly be upgraded to a homotopy
equivalence.

\noindent{\bf\S9.3 Identifying Components:}~The equivalence in 9.5
gives a homology equivalence at the level of components. We identify
these components for both $Par (M-N)$ and the space of sections.  For
the sake of simplicity we confine ourselves to the case $M-N$
parallelizable.

\noindent{\bf Lemma} 9.14:~{\sl Let $X$ be a connected topological space,
$M, N$ as above, $N\neq\emptyset$. Then
all components of $\map{}(M,N, \party (S^n,*))$
are homotopy equivalent.}

\noindent{\sc Proof}: 
Let $*\in N\subset M$ and identify $*$ with $N$ in $M/N$. 
Pick a map $f\in\map{}$ and observe that for a small disc
$D\in M-N$, $f(\partial D)$ is null homotopic in $\party (S^n,*)$
(since the latter is $n-1$ connected).
So $f_{|\partial D}$ extends out to a map of a sphere $S^n$ and if
we denote by $\#$ the connected sum, we have a map
$$\map{}^*(M/N,\party (S^n,*))\lrar\map{}^*((M/N)\# S^n,\party (S^n,*))$$
which takes one component to the next (here of course $(M/N)\#
S^n\simeq M/N$).  This map is a homotopy equivalence for it can be
reverted by attaching another sphere with reverse orientation.  
\hfill\za

\noindent{\bf Proposition} 9.15:~{\sl Let $\overline{M-N}$ be the
closure of $M-N$ and let $p\in\partial
(\overline{M-N}))\neq\emptyset$. Then
$$S: \party(\overline{M-N},p)\fract{S}{\ra 3}\map{0}(M,N, \party (S^n,*))$$
is a homology equivalence.
Here $\map{0}$ stands for the component of null-homotopic maps.}

\noindent{\sc Proof:} 
Recall (\S7) that the construction of $Par$ involved a choice of 
multiconfigurations $\vec\zeta_i\in U_i\subset U$ where $U$ is an open
collar around $N$ and the $U_i$ form a nested sequence contracting to $U$.
The stabilization maps $+\vec\zeta_i$ are easily seen to commute with
scanning and hence we get a commuting diagram ($\map{c}$ and $\map{c'}$
are some components)
$$\matrix{
  \party (M-U_i)&\fract{S_i}{\lrar}&
  \map{c}(M/U_i,\party (S^n,*))\cr
  \decdnar{+\vec\zeta_i}&&\decdnar{}\cr
  \party (M-U_{i+1})&\fract{S_{i+1}}{\lrar}&
  \map{c'}(M/U_{i+1},\party (S^n,*)).
\cr}\leqno{9.16}$$
Observe that $\party (M-U_i)\cong\party (M-N)$ for all $i$ and up to
homotopy we have $\map{c}=\map{0}$. We then get the homotopy commutative
diagram
$$\matrix{
  \party (M-N)&\fract{S}{\lrar}&
  \map{0}(M,N,\party (S^n,*))\cr
  \decdnar{+\vec\zeta_i}&&\decdnar{}\cr
  \party (M-N)&\fract{S}{\lrar}&
  \map{0}(M,N,\party (S^n,*)).
\cr}$$
which yields in the limit a map of components
$$\party (\overline{M-N},p )\fract{S}{\ra 2}\map{0}(M,N, \party (S^n,*)$$
and this must be a homology equivalence.
\hfill\za

\noindent{\bf Remark} 9.17: Since $\party (S^n,*)$ is $n-1$ connected
(proposition 6.10) it follows that
$$\pi_0\map{}(M,N, \party (S^n,*)) = [M/N,\party (S^n,*)]_*$$ 
and hence that the connected components of the corresponding mapping
space (and consequently of $Par(M-N)$) are indexed by maps of
$H_n(M,N;\bbz )$ into $H_n(\party (S^n,*))$.


\vskip 10pt
\noindent{\bf\Large\S10 Duality on Manifolds}

As mentioned in the introduction, theorem 9.5 admits a strengthening
when $\party=\spy$. This last functor is a homotopy functor on the one
hand, and on the other it takes values in abelian monoids. We start
with some standard results.

First we point out that since $K(\bbz,n)$ has the homotopy type of an
abelian group, then so does the space of maps $\map{}(X,K(\bbz, n))$
and for connected $X$, all components of $\map{}(X, K(\bbz, n))$ are
homotopy equivalent.  A classical result (attributed to Moore) asserts
that any abelian topological group is a product of Eilenberg-MacLane
spaces. It remains to determine what these EM spaces are for the case
of $\map{}(X, K(\bbz, n))$ and this is exactly the content of the
following theorem of Thom (cf. [NS])

\noindent{\bf Theorem} 10.1: (Thom)~{\sl Let $X$ be connected, $\pi$
an abelian group and $n>0$. Then
$$\displaystyle\map{}(X,K(\pi,n))\simeq \prod_{0\leq i\leq n}
K(H^{n-i}(X,\pi ),i)$$
and each component is given by the sub-product $1\leq i\leq n$ in the
expression above.}

Let $X=K(\bbz, n)$ and consider the subspace $Aut(K(\bbz,n))\subset
\map{}(K(\bbz, n), K(\bbz, n))$ of self-homotopy equivalences of
$K(\bbz,n)$. This is an abelian subgroup and hence is also a product
of EM spaces. Our next proposition is an earlier result of May [Ma] 
which we state and prove in ``non-simplicial'' terms.

\noindent{\bf Proposition} 10.2:~{\sl We have the following 
commutative diagram of inclusions and equivalences
$$\matrix{K(\bbz, n)\times Aut(\bbz)&\hookrightarrow&K(\bbz, n)\times
Hom(\bbz, \bbz)\cr
\decdnar{\simeq}&&\decdnar{\simeq}\cr
Aut(K(\bbz,n))&\hookrightarrow&\map{}(K(\bbz, n), K(\bbz, n))\cr}$$}

\noindent{\sc Proof}: To simplify notation we write $K_n:=K(\bbz, n)$.
From 10.1 and since $H^{n-i}(K_n;\bbz )=\bbz$ when $i=0$ and
zero otherwise, we get
$$\map{}(K_n, K_n)\simeq K(H^0(K_n;\bbz ), n)\times
K(H^n(K_n;\bbz ),0) \simeq K_n\times Hom(\bbz, \bbz )$$
(here of course $H^n(K_n;\bbz) =$ Hom$(H_n(K_n;\bbz ),\bbz )$
= Hom$(\bbz, \bbz)\cong \bbz$.)
The equivalence above can be explicitly contructed as follows.  We
pointed out earlier that $K(\bbz, n) \simeq \spy (S^n,*)$ (this
equivalence can be seen in many ways; cf. [DT] or [M]) and the abelian
monoid structure on $K_n=K(\bbz,n)$ is induced from the symmetric product
pairing (which we write additively). Given a map $f:\bbz\lrar\bbz$
determined by an integer $k$, we can consider the $k$-fold map
$S^n\lrar S^n$ and extend it out (additively) to a map $(k):\spy
(S^n,*)\rightarrow\spy (S^n,*)$ and hence to an element $(k)\in
\map{}(K_n, K_n)$. On the other hand, $K_n$ maps
to the translation elements in $\map{}(K_n,K_n)$ and the
product map $(x, k)\mapsto T_x+ (k)$ induces the equivalence
$K_n\times Hom(\bbz, \bbz)\rightarrow\map{}(K_n, K_n)$.
The homotopy inverse sends $f\in\map{}(K_n, K_n)$ to $(f(x_0), \deg f)$
where $x_0\in K_n$ is the basepoint and $\deg f$ is the degree of the
induced map at the level of $\pi_n$.

Note at this point that since $H_n(S^n)\cong\pi_n(\spy (S^n,*))$, the
map $(k)$ induces multiplication by $k$ at the level of $\pi_n$ and so
$(k)$ is a homotopy equivalence if and only if $k=\pm 1$, in which
case multiplication by $k$ is in $Aut(\bbz )$. Notice also that an
element in $K_n$ acting by translation can be homotoped to the
identity and hence the map $T: K_n\rightarrow\map{}(K_n,K_n)$ factors
through $Aut(K_n)$.  These two facts put together show that the diagram in
10.2 commutes. It remains to show that the left vertical map is an
equivalence but it is not hard to see that the right-hand
equivalence we just described restricts to $Aut(K_n)$
and the proposition follows. \hfill\za

\noindent{\bf Remark} 10.3: We can replace $\bbz$ by any abelian group
$G$ in 10.2 above and prove similarly that $Aut(K(G,n))\simeq
K(G,n)\times Aut (G)$. At the level of simplicial groups,
$Aut(K(G,n))$ is given as a semi-direct product of $Aut(G)$ and
$K(G,n)$ (May). When $G=\bbz$, $Aut(\bbz )\cong\bbz_2$ and $Aut(K_n)$
consists of two copies of $K_n$ (consisting resp. of ``orientation''
preserving and orientation reversing homotopy equivalences).

\noindent{\bf Theorem} 10.4:~{\sl Let $M\in\cat_n$. Then the bundle
$K(\bbz,n)\rightarrow E_{\spy}\rightarrow M$ is trivial if and only if
$M$ is oriented.}

\noindent{\sc Proof:} The bundle $E_{\spy}$ is classified by a map
$M\lrar BAut(K(\bbz,n))$ and at the level of spaces we get a (trivial)
fibration
$$K(\bbz,n+1)\lrar B\left(Aut (K(\bbz,n))\right)\lrar B\left(Aut
(\bbz)\right).$$
The classifying map $f:M\lrar BAut(K(\bbz,n))$ lifts to $K(\bbz,n+1)$
if and only if the composite $M\rightarrow B(Aut(\bbz))$ is null
homotopic or equivalently if the induced map $\phi: \pi_1(M)\lrar Aut
(\bbz)$ is trivial.  The action of $\pi_1(M)$ on $\bbz$ described by
the map $\phi$ corresponds to the action of $\pi_1(M)$ on $\bbz
=\pi_n(K(\bbz, n))$ in the bundle in 10.4 (this follows directly from
the many facts stated in the proof of 10.2).  But $M$ being oriented,
the tangent bundle $\tau M$ (and hence its compactified counterpart
$\widehat{\tau M}$) is trivial over the 1-skeleton.  Consequently,
$E_{\spy}$ restricted to the one skeleton of $M$ is also trivial and
so is the action of $\pi_1(M)$ on the fiber.  Namely, $\pi_1(M)$ acts
trivially on $\pi_n(K(\bbz,n))=\bbz$ and as indicated above the map
$f$ must lift to a map $\tilde f: M\lrar K(\bbz,n+1)$. Since $M$ is
$n$ dimensional, $\tilde f$ is null-homotopic and $E_{\spy}$ is
trivial.

To prove the other easier direction, suppose $E_{\spy}$ is trivial
that is $E_{\spy}\simeq K(\bbz,n)\times M$. The inclusion
$\widehat{\tau M}\subset E_{\spy}$ composed with projection
yields a map of $\widehat{\tau M}\rightarrow K(\bbz,n)$ and hence
a Thom class in $H^n(\widehat{\tau M};\bbz )$.
This is equivalent to giving an orientation class for
$M$ and the proposition follows.  \hfill\za

\noindent{\bf Theorem} 10.5:~{\sl Let $N\hookrightarrow M$ be a closed
ANR of a closed, oriented manifold $M\in\cat_n, n\geq 2$. Then
$$\spy (M-N,*)\fract{\simeq}{\ra 2}\map{0}(M,N, \spy (S^n,*)).$$}

\noindent{\sc Proof:} 
Here of course and since $E_{\spy (M-N)}$ is trivial, the space
of sections and the space of maps into the fiber coincide.
The homotopy equivalence is a consequence of 1.3 (or 9.13.)
\hfill\za

\noindent{\bf Remark} 10.6: When $M$ is parallelizable, the map $S$ has
the following alternate description.  Start with $M$ compact and for
each $x\in M$ choose an open ball $D_x\subset M$ containing $x$ and
such that $D_x/\partial D_x\simeq S^n$ canonically. The quotient maps
$M\lrar D^x/\partial D_x = M/M-D_x\simeq S^n$ give rise to maps
$$s_x:M\lrar S^n\hookrightarrow\spy (S^n,*),~\forall x\in M$$
and hence to a correspondence
$s: M\lrar\map{}(M,\spy (S^n,*))$ which extends additively to
$${\bar s}: \spy (M)\lrar\map{}(M, \spy (S^n,*)).$$
It isn't hard to see that ${\bar s}\simeq S$ (Another variation
on this construction is given in \S12.)

A direct consequence of proposition 2.7 and from the fact that
$\spy (-)$ is a homotopy functor, it follows that $\pi_*(\spy (-))$
defines a homology theory and a well-known theorem of Dold and Thom
identifies it with ordinary singular homology theory; i.e.
$$\spy (X,*)=\prod_iK({\tilde H}_i,(X;\bbz ),i)\leqno{10.7}$$
Combining 10.1 with 10.7 we get the equivalence
$$\prod_iK({\tilde H}_i,(M-N;\bbz ),i)\simeq \prod_{1\leq i\leq n}
K(H^{n-i}(M/N,\bbz ),i)$$
from which we easily deduce our main application

\noindent{\bf Corollary} 10.8: (Alexander-Poincar\'e Duality)~{\sl Let
$N\hookrightarrow M$ be a closed ANR in an orientable manifold $M$ of
dimension $n$. Then ${\tilde H}_i(M-N;\bbz )\cong H^{n-i}(M,N;\bbz
).$}

Similarly, considering the equivalence $\spy
(M,*)\simeq\map{c}(M,\partial M, K(\bbz,n))$ for $M$ compact with
boundary yields

\noindent{\bf Corollary} 10.9 (Lefshetz-Poincar\'e Duality)~{\sl Let
$M$ be compact with boundary, of dimension $n$, and suppose $int M$ is
orientable. Then $H_q(M)\cong H^{n-q}(M,\partial M)$.}

\noindent{\bf Example} 10.10: The classical Alexander duality is stated
as follows. Let $X$ be a finite complex embedded in $S^n$ ($n\geq 1$).
By 10.8 we have that
$H^{n-i}(S^n,X)\cong {\tilde H}_i(S^n-X)$ and 
the relative sequence for the pair $(X,S^n)$ shows that 
$${\tilde H}_i(S^n - X)\cong H^{n-1-i}(X;\bbz ).$$
When $i$ corresponds to one less the ``codimension'' of $X$ in $S^n$,
the isomorphism above has a very nice geometric interpretation
(see [KT] for example). Suppose $X=M$ is a smooth closed $m$ manifold
embedded in $S^n$. It has a unit normal sphere bundle
$S^{n-m-1}\rightarrow \nu (M)\rightarrow M$ and the homology class of
the fiber in ${\tilde H}_{n-m-i}(S^n - X)$ is dual to the cohomology
orientation class in $H^m(M;\bbz )$.  Note in this case that the class
in ${\tilde H}_{n-m-i}(S^n - X)$ is spherical.



\vskip 10pt
\noindent{\bf\Large\S11 Applications}

\noindent{\bf\S11.1 On Theorems of McDuff and Segal}:~
As pointed out in the introduction, the configuration space functor
$\config$ has been studied in [S2] and [McD1] where special versions
of theorem 1.1 have been proved. In this subsection, we extend their
results in several directions.

Consider the subspace of $C^{(k)}(M)\subset C(M)^k$ consisting of
tuples of configurations which are pairwise disjoint. More explicitly
$$C^{(k)}(M) = \{(\zeta_1,\ldots, \zeta_k)\in C(M)^k~|~\zeta_i\cap\zeta_j
=\emptyset, i\neq j\}.$$
It is direct to see that $C^{(k)}(M)$ is a particle space 
and hence for parallelizable pairs $(M,N)$ we have (theorem 9.15)  
$$H_*(C^{(k)}(M-N);\bbz )\fract{S_*}{\ra 2}
H_*(\map{}(M,N, C^{(k)}(S^n,*));\bbz )$$

\noindent{\bf Lemma} 11.1:~{\sl Let $\bigvee^{k}S^n$ denote 
the $k$-th wedge, $n\geq 1$. Then $C^{(k)}(S^n,*)\simeq \bigvee^{k}S^n$.}

\noindent{\sc Proof:} As in [S1], we let $C^{(k)}_{\gre}(S^n,*)$ be
the open set of $C^{(k)}(S^n,*)$ consisting of multiconfigurations
$(\zeta_1,\ldots, \zeta_k)$ such that at least $k-1$ such particles
are disjoint from the closed disk $U_{\epsilon}$ of radius $\gre>0$
about the south pole $*$.  Notice that there is a radial homotopy,
injective on the interior of $U_{\epsilon}$) that expands the north
cap $U_{\epsilon}$ over the sphere and takes $\partial U_{\epsilon}$
to $*$. Such an expansion retracts $C^{(k)}_{\gre}(S^n,*)$ to the
wedge product $C(S^n,*)\vee\cdots\vee C(S^n,*)$. Now since
$C^{(k)}(S^n,*)$ is the union of the $C^{(k)}_{\gre}(S^n,*)$ for
$\gre>0$, we get that $C^{(k)}(S^n,*)\simeq\bigvee^kC(S^n,*)$.

It remains to show that $C(S^n,*)\simeq S^n$. Here too we consider the 
subspace 
$$C_{\gre}(S^n,*)=\left\{ D\in C(S^n,*)~|~D\cap U_{\epsilon}=\{\hbox{at most
one point}\}\right\}$$
where $U_{\gre}$ is an epsilon neighborhood of the north pole (again
the south pole corresponds to $*$). Then radial expansion of
$U_{\gre}$ ($N$ is fixed) maps $(U_{\gre},\partial U_{\gre})$ to
$(S^n, *)$ (and is injective on $U_{\gre}$ hence extending to $C$).
The one point configurations in $U_{\gre}$ now produce a
homeomorphism $C_{\gre}(S^n,*)\simeq S^n$ and since again
$C(S^n,*)=\bigcup_{\gre} C_{\gre}(S^n,*)$ the lemma follows.
\hfill\za

\noindent{\bf Proposition} 11.2:~{\sl Let $M\in\cat_n$ be a closed
manifold and $N\subset M$ such that $(M,N)$ is parallelizable. Then
$$S_*:H_*(C^{(k)}(M-N))\cong H_*(\map{}(M,N, \bigvee^kS^n)).$$}

When $M-N=\bbr^n\cong D^n$, $D^n$ here is the closed unit disc, then
components of $C^{(k)}(\bbr^n)$ can be identified with the direct limit
$\config(D^n,p)$ constructed in 7.5. We have

\noindent{\bf Corollary} 11.3:~{\sl The scanning map $S:
C^{(k)}(\bbr^n)\lrar \Omega^n(\bigvee^kS^n)$ induces a homology
isomorphism. When $k=1$ we recover the following classical result of
Segal [S2]
$$
\displaystyle E_*:H_*(\config(D^n,p);\bbz )\cong
H_*\left(\lim_{\rightarrow\atop i}C_i(\bbr^n);\bbz \right)
\fract{\cong}{\ra 2} H_*(\loop^n_0S^n;\bbz ).
$$}

\noindent{\bf Example} 11.4: It can be checked (exactly as in 11.1) that
$DDiv^k(S^n,*)\simeq\bigvee^kK(\bbz,n)$ and that the following commutes
(up to homotopy)
$$\matrix{
   C^{(k)}(M-N)&\fract{S}{\lrar}&\map{}(M, N,\bigvee^kS^n)\cr
   \decdnar{\subset}&&\decdnar{}\cr
   DDiv^k(M-N)&\fract{\simeq}{\lrar}&\map{}(M, N,\bigvee^kK(\bbz,n))
\cr}$$
where $M$ and $N$ are as in the statement of theorem 1.1. 
We quickly remind the reader that $DDiv^k(M)$ is the set of
$k$-tuples of positive divisors which are pairwise disjoint. We finally
point out that the right
vertical map in the diagram is induced from the inclusion
$S^n\hookrightarrow K(\bbz, n)$ and the homotopy equivalence at the
bottom follows from the fact that $\pi_1(DDiv^k(\bbr^n))$ is abelian
(which is left for check to the reader).\hfill\za

\noindent{\bf Example} 11.5:
{\sl (Spaces of positive and negative particles)}~
[McD1] also introduces the functor $C^{+\over}$ discussed in \S1. 
This is given as the quotient of $C(M)\times C(M)$ with the
relation 
$$ (\langle x, x_1,\ldots, x_n\rangle,\langle x,y_1,\ldots
,y_m\rangle ) \sim_R (\langle x_1,\ldots, x_n\rangle,\langle
y_1,\ldots ,y_m\rangle ).  
$$
One can show that $C^{+\over}$ abelianizes fundamental group
and hence theorem 1.3
applies.  Observe that since $C(S^n,*)\simeq S^n$, it follows that
$C^{+\over}(S^n,*) \simeq S^n\times S^n/\Delta$ where $\Delta (S^n)$
is the diagonal copy of $S^n$ in $S^n\times S^n$. The following 
homotopy equivalence is a
special case of [McD1] or of theorem 1.3
$$C^{+\over}(\bbr^n)\simeq\Omega^n\left(S^n\times S^n/\Delta (S^n)\right).$$

\noindent{\bf\S 11.2 Symmetric products with bounded multiplicities}:~
In this subsection we prove theorem 1.5 in the introduction. Recall
that $\spy_d$ was defined as the particle functor of the first kind
$$\spy_d(M) =\{\sum n_ix_i\in\spy (M)~|~n_i\leq d\}.$$
We first need the following analog of 11.1.

\noindent{\bf Lemma} 11.6:~ {\sl There is a homotopy equivalence
$\spy_d(S^k,*)\simeq\sp{d}(S^k)$.}

\noindent{\sc Proof:} Let $*\in S^k$ and $U_{\epsilon}$ 
be as in 11.1, and let
$W_{\epsilon}$ be the subspace consisting of
$\langle x_1,x_2,\ldots, x_n\rangle\in \spy_d(S^k,*)$ such that at most
$d$ points in the tuple lie inside $U_{\epsilon}$. 
By definition of $\spy_d(S^k,*)$ each of its elements must fall
into a $W_{\epsilon}$ for some $\epsilon$ and hence
$$\spy_d(S^k,*)\simeq \bigcup_{\epsilon}W_{\epsilon}.$$
Now using the radial retraction of 11.1, it is clear that each
$W_{\epsilon}$ is homotopically $\sp{d}(S^k)$ (since by
taking a configuration $\langle x_1,x_2,\ldots, x_n\rangle$ and
shrinking (at least) $n-d$ points of it to basepoint $*$, we end up in
$\sp{d}(S^k)$.) The lemma follows.
\hfill\za

\noindent{\bf Theorem} 11.7:~{\sl Let $M$ and $N$ be as in 1.1. Then 
$$S: SP_d(M-N)\lrar\map{}(M,N\cup\partial M,\sp{d}(S^n))$$
is a homotopy equivalence whenever $d>1$ and a homology equivalence
when $d=1$.}

\noindent{\sc Proof}: Let $X=M-N$. The claim
amounts to showing that $\pi_1(\sp{n}_d(X))$
is abelian when $n>1$ and $d>1$.
We know already (9.13) that $\pi_1(\sp{n}(X))$ is abelian for $n>1$.
Since $H_1(\sp{n}(X);\bbz)\cong H_1(\sp{n+1}(X);\bbz )$, it follows
that the inclusion $\sp{2}(X)\hookrightarrow\sp{n}(X)$ for $n\geq 2$
induces an isomorphism in fundamental group. Consider at this point
the commutative diagram
$$\matrix{
\sp{2}_d(X)&\hookrightarrow&\sp{2}(X)\cr
\decdnar{\subset}&&\decdnar{\subset}\cr
\sp{n}_d(X)&\hookrightarrow&\sp{n}(X)\cr}$$
Any element $\alpha\in\pi_1(\sp{n}_d(X))$ factors through
the subset $\sp{2}_d(X))$ in $\sp{2}(X)$. But for $d>1$, 
these last two spaces coincide and since $\pi_1(\sp{2}(X))$ is 
abelian, the claim follows.
\hfill\za

\noindent{\sc Corollary} 11.8: Restricting to the case $M=D^n$ the closed
unit disc, $N=\emptyset$, we recover 1.5 in the introduction (see also
5.2). When $X=\bbc\cong\bbr^2$, the space $\sp{n}_d(\bbc )$ can
be identified with the space of monic polynomials $p$ 
of degree $n$, $p(z)=(z-x_1)\cdots (z-x_n)$ such that $p$ has no roots
of multiplicity greater than $d$. Notice in this case that
$\sp{d}(S^2)$ is diffeomorphic to the $d$th complex projective space
$\bbp^d$ and hence we obtain the following corollary also proved in [GKY]

\noindent{\bf Corollary} 11.9 [GKY]:~{\sl There is a correspondence
$$\left\{{\hbox{Monic complex polynomials of degree $n$}
\atop\hbox{and roots of multiplicity $d>1$}}\right\}\fract{}{\lrar}
\Omega^2_0\bbp^d$$
which is a homotopy equivalence in the direct limit when $n\lrar\infty$.}

\noindent{\bf Remark} 11.10: One can prove more in this case (as [GKY]
do) by showing that the correspondence above is a homotopy equivalence
through a range. This is a good place to point out that our main
theorem 1.1 is quite likely to have an unstable version which would
state that scanning $S$ is a homology equivalence through a range
increasing with the multidegree of the $Par$ spaces.

\noindent{\bf\S11.3 Rational curves on toric varieties and a
theorem of Guest}:~
A toric variety $V$ is a projective variety that can be defined by
equations of the form ``monomial in $z_0,\ldots, z_n$= monomial in
$z_0,\ldots, z_n$". As an example, consider the quartic
$$M_2=\{ [z_0:z_1:z_2,z_3]\in\bbp^3~|~z_2^2=z_1z_3\}.$$

A rational curve on $V$ is a holomorphic image of $\bbp^1=S^2$ in $V$
and one is interested in studying the space of all such curves. The
interest here stems from the relevance of these spaces to problems in
Gauge theory, Sigma models for physicists and even Control theory for
engineers (cf. [BHMM], [C$^2$M$^2$] and [S1] for a general discussion
of the subject).

We denote by $\hol{}(\bbp^1, V)$ the space of all holomorphic maps
from $\bbp^1$ into $V$.  As is customary, the study of this space
proceeds  by first restricting attention to the
subspace of based maps (which consists of maps that fix a given
basepoint.) Choosing $x_0\in\bbp^1$ and $*\in V$, we let
$\hol{}^*(\bbp^1, V)$ be the subspace of $f:\bbp^1\rightarrow V$
such that $f(x_0)=*$. It has to be pointed out that the topology of
$\hol{}^*(\bbp^1,V)$ could vary with the choice of the basepoint $*$
(unless for example $V$ is homogeneous).

It turns out that for a generic choice of a basepoint $*\in V$, a map
$f\in\hol{}^*(\bbp^1,V)$ admits a representation by polynomials. More
precisely, given $f:\bbp^1\lrar V$ holomorphic, the composite
$$\bbp^1\fract{f}{\lrar}V\hookrightarrow \bbp^n~~~\hbox{(for some}~n)$$
is also holomorphic and so $f$ can be represented by the map
$[p_0(z):\ldots: p_n(z)]$ where the $p_i(z)$ satisfy the same set of
equations as $V$ and of course have no roots in common. Notice also that
when $f$ is basepoint preserving, the $p_i$ can be chosen to be monic
(and hence are uniquely determined).

\noindent{\bf Example} 11.11 (Guest): Consider the quadric curve $M_2$
described earlier. It can be seen that $M_2$ is smooth but at the one
point $[1:0:0:0]$.  A (based) rational curve $f: S^2\rightarrow M_2$,
sending the north pole to any point other than this singular point,
has a representation in terms of a $4$-tuple of polynomials
$(q_1,q_2,q_3,q_4)$ that are coprime, monic and satisfying the equation
$q_3^2=q_2q_4$. The map $f$ is therefore equivalent to the choice of
four monic polynomials $p_1,p_2,p_3$ and $p_4$ such that
\begin{eqnarray*}
   &&(q_1,q_2,q_3,q_4)= (p_4, p_1^2p_2 , p_1p_2p_3 ,p_2p_3^2)\\
   &&(p_1,p_3)=1, (p_2,p_4)=1\\
   &&\deg p_4 =\deg p_1^2p_2 = \deg p_1p_2p_3 = \deg p_2p_3^2 =d
\end{eqnarray*}
where $d$ is the degree of $f$. This last representation can be
reformulated in terms of divisors $D_1,D_2,D_3$ and $D_4$
given by the roots of the $q_i$ and hence satisfying
\begin{eqnarray*}
  && D_1\cap D_2\cap D_3\cap D_4=\emptyset, ~\deg D_i=d\\
  && D_2, D_3, D_4 ~\hbox{are of the form}~ \eta + 2\phi, \eta + \phi +
  \chi, \eta + 2\chi 
  ~\hbox{where}~ \phi\cap\chi =\emptyset.
\end{eqnarray*}

Generally, the $p_i$'s one associates to $f\in\hol{}^*(\bbp^1,V)$ 
being monic, their root data totally
determine the map $f$. For a general toric variety $V$, which we
assume to be non-singular (the singular case is a little more
intricate but can still be treated analgously), a rational map
$f:\bbp^1\lrar V$ will have a multidegree $D$ associated to it where
$$D = (d_1,\ldots, d_p)\in \pi_2(V)\cong\bigoplus_{i=1}^p\bbz$$
and this multidegree parametrizes components of $\hol{}^*(\bbp^1,V)$.
We say $D\rightarrow\infty$ if all the components $d_i$ tend to
infinity.

\noindent{\bf Lemma} 11.12:~{\sl There is a homeomorphism
$\hol{D}^*(S^2,V)\cong Par_D (S^2-{\infty})$ for some particle space
$\party (S^2-*)$, sending $f\in\hol{d}^*(S^2,V)$ to the roots of the
$p_i(z), 0\leq i\leq n$ in its polynomial representation.}

\noindent{\sc Proof:} The proof is direct since if two polynomial
representations given by $p_i$ and $p_i'$, $1\leq i\leq n$, have root
data lying in disjoint sets, then their products $p_ip_i'$ will give
rise to another representation describing a new holomorphic map
$S^2\rightarrow V$.  \hfill\za

We can up to homeomorphism contruct stabilization maps
$$\hol{D}^*(S^2,V)\lrar \hol{D+D'}^*(S^2,V)\leqno{11.13}$$ 
as in \S7.  This induces stabilization maps at the level of
$\party_D(S^2-\infty)$ and the direct limit is a component of
$Par(S^2-\infty )$ (see \S5).

\noindent{\bf Theorem} 11.14:~(Guest)~{\sl Let $X$ be a projective
toric variety (non-singular). The inclusions
$i_D:\hol{D}^*(S^2,V)\lrar\loop^2_DV$ induce a homotopy equivalence 
when $D$ goes to $\infty$; i.e.  
$$\lim_{D\rightarrow\infty}
\hol{D}^*(S^2,V)\fract{{\lim\atop \rightarrow}i_D}{\ra 4}\loop^2_0V$$
where $\loop^2_0V$ is any component.}

\noindent{\sc Proof:} Arguments of Segal and Guest show that in this
general case scanning and the inclusion $i$ fits in a homotopy
commutative diagram
$$\matrix{
  \hol{D}^*(S^2,V)&\fract{i_D}{\ra 2}&\map{D}^*(S^2,V)\cr
  \decdnar{}&&\decdnar{\simeq}\cr
  Par_c(S^2-\infty)&\fract{S}{\ra 2}&\map{c}^*(S^2,\party (S^2,*)).
\cr}$$
where from above the map $\hol{D}^*(S^2,V )\rightarrow
Par_c(S^2-\infty )$ can be identified with the map of $\hol{D}^*$ into
the direct limit of the system in 11.13 (note that $\map{c}^*$ denotes
any component of $\map{}^*(S^2,\party (S^2,*))$ and they're all
homotopy equivalent by 9.14).  The scanning map $S$ at the bottom will
be a homotopy equivalence according to 1.3 if we can show that
$\pi_1(Par_c(S^2-\infty))$ is abelian.  It is shown in ([BHMM],
corollary 9.9) that $\pi_1(\hol{D}^*(S^2,V))$ is abelian for $D$
consisting of multidegrees $(d_1,\ldots, d_p)$ with $d_i\geq 2$.
Moreover for $D$ and $D'$ with this property,
$\pi_1(\hol{D}^*(S^2,V))\cong\pi_1(\hol{D'}^*(S^2,V))$ hence implying
that in the direct limit $\pi_1(\hol{}^*(S^2,V))$ is well defined and
abelian. The claim now follows.  \hfill\za


\vskip 10pt
\noindent{\bf\Large\S12 Spanier-Whitehead Duality}

The ideas of the previous sections can be adapted to prove
Spanier-Whitehead duality for general homology theories $h_*$ and for
any finite type CW complex $X$. The material below is known in some
form or another and we include it in this section for completeness.

As a start we denote by $CW$ the category of connected finite CW
complexes.  For a given $X\in CW$, we let $D(X,k)$ be its
Spanier-Whitehead dual (or $S$-dual). An $S$-dual always comes
equiped with a map $X\wedge D(X,k)\lrar S^k$ (see [CM]).

One can construct the $S$-dual of any $X\in CW$ very
concretely. Indeed, since $X$ is finite, it embeds in some big enough
sphere $S^k$. The complement $Y=S^k-X$ can now be chosen as a
spanier-whitehead dual for $\Sigma X$; i.e. $Y=D(\Sigma X,k)$
([CM]). The $S$-dual of $X$ can then be taken to be $\Sigma D(\Sigma
X,k)$.  It follows for instance that the $S$ dual of $S^n$ is
$S^{k-n}$.

Given a connective $\Omega$ spectrum ${\bf E} = \{E_i, i=1,\ldots\}$, 
we have that 
$$E_0=\lim_m\Omega^mE_{m} \equiv \Omega^{\infty}{\bf E}$$
and more generally $E_n=\Omega^{\infty}(S^n\wedge {\bf E})$.
The functor $\Omega^{\infty}$ is a functor from spectra to spaces.

We can associate to ${\bf E}$ the functor $F_{\bf E}$ defined as follows
$$F_{\bf E}: X\mapsto F(X)=\Omega^{\infty}({\bf E}\wedge X).$$
Notice that by definition $F_{\bf E}(S^n)=E_n$
Notice also that
\begin{eqnarray*}
    \pi_i(F_{\bf E}(X))&=&[S^i, \Omega^{\infty}(X\wedge {\bf E})]
    = \lim_n[S^i,\Omega^n(E_n\wedge X)] \\
    &=&\lim_n\pi_{i+n}(E_n\wedge X)=h_i(X)
\end{eqnarray*}
where $h_*$ is the generalized homology theory associated to $\bf E$.

\noindent{\bf Remark} 12.1: Generally, given a spectrum ${\bf E}$, we
denote by $\Omega^{\infty}$ the functor obtained as the composite of
the functor which converts any spectrum into an equivalent $\Omega$
spectrum ${\bf E}$ followed by the functor which passes from ${\bf E}$
into the space $E_0$ (see [Ad],p:22).  Notice that $E_0$ doesn't
generally correspond to $E_0'$ (as the sphere spectrum
$\bf S^0$ does illustrate already). 

\noindent{\bf Theorem} 12.2:~{\sl Let $F=F_{\bf E}$ for some spectrum
$\bf E$. Then $\forall X\in CW$, there is a homotopy equivalence
$$S:F(X)\fract{\simeq}{\ra 2}\map{}^*(D(X,k), F(S^k)).$$}

\noindent{\sc Proof:} Let $X$ be a finite CW complex. Then $X\subset
S^k$ for some $k$ and $D(\Sigma X,k)=S^k-\Sigma X$. Since $X$ and $D(\Sigma
X,k)$ are disjoint, we can consider the map
$${\hat S}: X\times D(\Sigma X,k)\lrar S^{k-1},~(x,y)\mapsto {x-y\over |x-y|}\in S^k.$$
We can assume $X$ to be embedded in the positive quadrant in
$\bbr^n\subset S^k$ with the point at infinity $\infty\in S^k$
adjoined. This means that ${\hat S}(\infty, y)=1,\forall y\in D(\Sigma
X,k)$.  On the other hand and since $X$ is compact, it lies in a ball
$B\in S^k$. Choose a point $p\in D(\Sigma X,k)$ which is not in
$B$. The map ${\hat S}_{|X\times p}$ extends to $B\times p$ and since
$B$ is contractible we get an extension
$${\hat S}:X\times D(\Sigma X,k)\cup c(X\times p)\lrar S^{k-1}$$
where $c$ denotes the cone construction. 
It then follows that up to homotopy, the map ${\hat S}$ gives rise to
the map
$$X\wedge D(\Sigma X,k)\lrar S^{k-1}.$$ 
Suspending both sides yields a map
$X\wedge D(X,k)\lrar S^k$ and hence by adjointing a map
$${\hat S}: X\lrar\map{}^*(D(X,k),S^k)\leqno{12.3}$$
where the mapping space on the right is pointed. Of course we can
compose with the map
$i:\map{}^*(D(X,k), S^k)\rightarrow \map{}^*(D(X,k),F(S^k))$ induced
from the ``identity" $S^k\rightarrow F(S^k)$. Since
$F$ is an infinite loop functor, 12.3 composed with $i$ extends to
the desired map
$$S: F(X)\lrar \map{}^*(D(X,k), F(S^k)).$$
We show that $S$ is a homotopy equivalence by inducting on cells of $X$.
Let $X^{(i)}$ denote the $i$-th skeleton of $X$ and consider the standard
cofibration $X^{(i-1)}\hookrightarrow X^{(i)}\lrar\bigvee S^i$.
Applying $\map{}^*(-,S^k)$ yields a fibration sequence and a homotopy
commutative diagram
$$\matrix{
  \prod F(S^{i-1})&\fract{\simeq}{\ra 2}&\prod \Omega^{k-i+1}F(S^k)\cr
  \downarrow&&\downarrow\cr
  F(X^{(i-1)})&\fract{S}{\ra 2}&\map{}^*(D(X^{(i-1)},k), F(S^k))\cr
  \downarrow&&\downarrow\cr
  F(X^{(i)})&\fract{S}{\ra 2}&\map{}^*(D(X^{(i)},k), F(S^k)) 
\cr}$$ 
The left hand vertical sequence is a quasifibration since $F$ is a
homology theory.  The top horizontal map is an equivalence since
$\Omega F(S^l)\simeq F(S^{l-1})$ while the bottom map is an
equivalence by induction. This then implies that the middle map $S$ is
also an equivalence and the proof follows.  \hfill\za

\noindent{\sc Example} 12.4:~By Dold-Thom, we know that
$\spy (-)$ is associated to the Eilenberg-MacLane spectrum $K(\bbz )$ 
(i.e. $\pi_*(\spy (-))\cong H_*(-;\bbz )$), while the
functor $Q(-)$ given by
$$QX=\Omega^{\infty}\Sigma^{\infty}(X)$$
is known to be associated to the sphere spectrum 
(i.e. $\pi_*(Q(X))=\pi_*^S(X)$).
One has then the following homotopy equivalence (described in [C])
$$Q(X)\fract{\simeq}{\ra 2} \map{}^*(D(X,k), QS^k).$$

At this point, consider $A,B\in S^k$ for some large $k$. Recall that
$A$ is $n$ dual to $B$ if $A\cap B=\emptyset$ and each is a strong
deformation retract of the complement of the other.

\noindent{\bf Corollary} 12.5: (Spanier-Whitehead duality)
~{\sl Let $\bf E$ be a connective spectrum and let $h$ be the homology 
theory defined by $\bf E$; i.e. $h_*(X) = [S^0,E\wedge X].$
Suppose $A,B\in S^k$, $A$ and $B$ are $n$ dual. Then there is
an isomorphism
$$h_i(B)\cong h^{n-1-i}(A).$$}

\noindent{\sc Proof:} Let $\bf E$ be a connective spectrum with a
unit. We can choose $\bf E$ to be an $\Omega$ spectrum. Indeed if it
weren't such, then the spectrum representing the generalized homology
theory defined by $\bf E$ still is.  And so as far as homology is
involved, we could have chosen $\bf E$ to be an $\Omega$ spectrum to
start with.

Theorem 12.3 now shows that $F_{\bf E}(X)\simeq\map{}(D(X,k), F_{\bf
E}(S^k))$ and it follows that
\begin{eqnarray*}
   h_i(X)&=&\pi_i(F_{\bf E}(X)) = \pi_i\left(\map{}(D(X,k), 
   F_{\bf E}(S^k))\right)\\
   &=&[S^i\wedge D(X,k), F_{\bf E}(S^k)] = [D(X,k),\Omega^iE_k]\\
   &=&[D(X,k), E_{k-i}] = h^{k-i}(D(X,k)).
\end{eqnarray*}
Here we used the facts that
$h_i(X) = \pi_i(F_{\bf E}(X)$ and $F_{\bf E}(S^n)\simeq E_n$. 
This concludes the proof.
\hfill\za


\vskip 20pt
\addcontentsline{toc}{section}{Bibliography}
\bibliography{biblio}

\begin{thebibliography}{doC}

\bibitem{}[A] J.F. Adams, ``Infinite Loop Spaces", Annals of Math. Studies,
study {\bf 90}, Princeton.
\bibitem{} [BS] J.C. Becker, R.E. Shultz, ``Equivariant function spaces
and stable homotopy theory I'', Comm. Math. Helv. {\bf 49} (1974)
1--34.
\bibitem{}[B] C.F. Bodigheimer, ``Stable splittings of mapping spaces'', 
Algebraic topology, Proc. Seattle (1985), Springer lecture notes {\bf 1286}, 
174--187.
\bibitem{}[BCM] C.F. Bodigheimer, F.R. Cohen, R.J. Milgram, ``Truncated 
symmetric products and configuration spaces'', Math. Zeit., {\bf 214} (1993), 
179--216.
\bibitem{}[BHMM] C.P. Boyer, J.C. Hurtubise, B.M. Mann, R.J. Milgram, ``The 
topology of rational maps into generalised flag manifolds'', Acta Matematica,
{\bf 173} (1994) 61--101. 
\bibitem{}[C] F.R. Cohen,
``Fibration and Product Decompositions in Nonstable Homotopy Theory", 
Hadbook of Algebraic Topology, Elsevier 1995.
\bibitem{}[C$^2$M$^2$] F.R. Cohen, R.L. Cohen, B.M. Mann, R.J. Milgram,
``The topology of rational functions and divisors of surfaces", Acta Math.,
{\bf 166}(1991), 163--221.
\bibitem{}[CM] G. Carlsson, R.J. Milgram, ``Stable homotopy and
iterated loop spaces'', Handbook of algebraic topology, Elsevier 1995.
\bibitem{}[DT] A. Dold, R. Thom, ``Quasifaserungen and unendliche symmetrische 
produkte'', Ann. Math. {\bf 67}(1958), 239--281.
\bibitem{}[G] P. Gajer, ``Poincar\'e duality and integral cycles'',
Compositio Math. {\bf 98}, (1995) 193--203.
\bibitem{}[GKY] M.A. Guest, A. Koslowski, K. Yamaguchi, ``The space of
polynomials with roots of bounded multiplicity'', math-AT/9807053.
\bibitem{}[Gu1] M. Guest, ``On the space of holomorphic maps from the 
Riemann sphere to the quadric cone'', Quart. J. Math. Oxford (2){\bf 45}, 
(1994), 57--75.
\bibitem{}[Gu2] M. Guest, ``The topology of the space of rational
curves on a toric variety'', Acta Math.,  {\bf 174}, (1995) 119--145.
\bibitem{}[K1] S. Kallel, ``Divisor spaces on punctured Riemann
surfaces'', Trans. Am. Math. Soc. {\bf 350} (1998), 135--164.
\bibitem{}[K2] S. Kallel, ``An interpolation between homology and stable
homotopy'', preprint April 98.
\bibitem{}[KM] S. Kallel, R.J. Milgram ``The geometry of spaces of holomorphic
maps from a Riemann surface into complex projective space'', 
J. Diff. Geometry {\bf 47} (1997) 321--375.
\bibitem{}[Ko] A. A. Kosinski, ``Differentiable Manifolds'', volume {\bf 138},
 pure and applied mathematics, academic press.
\bibitem{} [KT] V. Krushkal, P. Teichner, ``Alexander duality, gropes and
link homotopy'', MSRI preprint 1997-057.
\bibitem{} [Ma] P. May, ``Simplicial objects in algebraic topology'',
Van Nostrand (1967).
\bibitem{}[McD1] D. McDuff, ``Configuration spaces of positive and negative
 particles'', Topology {\bf 14}(1975), 91--107.
\bibitem{} [McD2] D. McDuff, ``Configuration spaces'', K-Theory and
Operator Algebras proc., Athens, Georgia (1975). Springer lecture
notes in Math. {\bf 575}, 88--95.
\bibitem{}[Mi] J. Milnor, ``Differentiable structures'', lecture notes, 1961.
\bibitem{}[M] R.J. Milgram, ``The bar construction and abelian H-spaces",
Ill. J. Math. {\bf 11} (1967), 242--250.
\bibitem{}[MS] D. McDuff, G. Segal, ``Homology fibrations and the group
completion theorem'', Invent. Math., {\bf 31} (1976), 279--284.
\bibitem{}[NS] D. Notbohm, L. Smith, ``Rational homotopy of the space of 
homotopy equivalences of a flag manifold'', 1990 Barcelona conference on 
Alg. Top., Springer lecture notes in Math.
\bibitem{}[S1] G. Segal, ``The topology of spaces of rational functions",
Acta. Math., {\bf 143}(1979), 39--72.
\bibitem{}[S2] G. Segal, ``configuration spaces and iterated loop spaces",
Invent. Math., {\bf 21}(1973), 213--221.
\bibitem{}[S3] G. Segal, ``K-Homology theory and algebraic K-theory'',
K-Theory and
Operator Algebras proc., Athens, Georgia (1975). Springer lecture
notes in Math. {\bf 575}, 113--127.
\bibitem{}[S4] G. Segal, ``Categories and homology theories'', Topology
{\bf 13} (1974) 293--312.

\end{thebibliography}
\bibliographystyle{plain}

\vskip 40pt

\flushleft{
Sadok Kallel\\
Dept. of Math., \#121-1884 Mathematics Road\linebreak
U. of British Columbia, Vancouver V6T 1Z2 \linebreak
{\sc Email:} skallel@math.ubc.ca}

\end{document}